\documentclass[a4paper,12pt]{article}

\usepackage{amscd}
\usepackage{amsfonts}
\usepackage{amsmath}
\usepackage{amssymb}
\usepackage{amsthm}
\usepackage{ascmac}
\usepackage{here}
\usepackage{makeidx}
\usepackage{mathrsfs}
\usepackage{slashbox}
\usepackage{txfonts}

\allowdisplaybreaks
\makeindex

\newcommand{\N}{\mathbb{N}}
\newcommand{\Q}{\mathbb{Q}}

\newcommand{\A}{\mathscr{A}}

\newcommand{\Fil}{\mathbb{F}}

\newcommand{\M}{\mathscr{M}}

\renewcommand{\P}{\mathscr{P}}

\newcommand{\V}{\mathscr{V}}

\newtheorem{thm}{Theorem}[section]
\newtheorem{dfn}[thm]{Definition}
\newtheorem{lmm}[thm]{Lemma}
\newtheorem{prp}[thm]{Proposition}
\newtheorem{crl}[thm]{Corollary}

\newtheorem{exm}[thm]{Example}

\def\ens#1{\{ #1 \}}
\def\Ens#1{\left\{ \ #1 \ \right\}}
\def\set#1#2{\{ \ #1 \mid #2 \ \}}
\def\Set#1#2{\left\{ \ #1 \ \middle|\ #2 \ \right\}}
\def\m#1{{\rm #1}}
\def\v#1{| #1 |}
\def\n#1{\| #1 \|}

\title{Iwasawa Theory for Locally Profinite Groups}
\author{The University of Tokyo, Tomoki Mihara}
\date{}

\begin{document}

\maketitle
\begin{abstract}
We establish duality theory of $p$-adic unitary Banach representations of locally profinite groups. This is an extension of Iwasawa theory for profinite groups by P.\ Schneider and J.\ Teitelbaum. We also establish a criterion for an irreducibility of a unitary representation by certain simpleness of the Iwasawa module dual to it. Through the duality, the continuous induction of a unitary Banach representation of a closed subgroup is interpreted as a certain induction of Iwasawa module. It gives an explicit description of the dual of a continuous parabolic induction on $\m{GL}_n(\Q_p)$.
\end{abstract}

\tableofcontents

\section{Introduction}
\label{Introduction}

Let $k$ be a local field, and $k^{\circ} \subset k$ the ring of integral elements. In \cite{ST}, P.\ Schneider and J.\ Teitelbaum established duality theory between Banach $k$-linear representations of a profinite group $G$ and Iwasawa modules over the Iwasawa algebra $k[[G]]$. We extend this result to the case when $G$ is a locally profinite group such as $\m{GL}_n(\Q_p)$. While every Banach $k$-linear representation of a compact group is unitarisable, i.e.\ admits an equivalent norm stable under the group action, the same does not hold in the case when the group is not compact. Therefore we deal only with unitarisable Banach representations.

\vspace{0.1in}
We recall basic notions of topological modules in \S \ref{Topological Module}. At the same time, we explain the convention in this paper. We define the Iwasawa algebra $k[[X]]$ of a topological space $X$ in \S \ref{Iwasawa Algebra}. Beware that $k[[X]]$ is not a ring in general. When $X$ is a locally profinite group $G$, we define a $k$-algebra structure on $k[[G]]$ associated to the group structure of $G$. This is a generalisation of the well-known Iwasawa algebra of a profinite group. We introduce the notion of bitopological $k^{\circ}$-algebra, and verify that a $k^{\circ}$-subalgebra $k^{\circ}[[G]] \subset k[[G]]$ admits a natural structure of a bitopological $k^{\circ}$-algebra.

\vspace{0.1in}
Suppose that $G$ is a locally profinite group. The duality theory of Banach $k$-linear representations of a profinite group and Iwasawa modules by P.\ Schneider and J.\ Teitelbaum is based on the duality between Banach $k$-vector spaces and compact Hausdorff flat linear topological $k^{\circ}$-modules. We also study the relation between a strongly continuous action of $G$ on a Banach $k$-vector space (resp.\ a compact Hausdorff flat linear topological $k^{\circ}$-module) and a $k^{\circ}$-algebra homomorphism from $k^{\circ}[[G]]$ to the ring of continuous automorphisms in \S \ref{Banach Spaces and Contractive Homomorphisms} (resp.\ \S \ref{Compact Linear Topological Modules}). Here we deal only with a unitary representation. As a consequence, we obtain duality between the categories $\m{Rep}_{G}(\m{BanC}^{\m{crt}}(k))$ and $\m{TMod}_{\m{flat,lin}}^{\m{cpt,sep}}(k^{\circ}[[G]])$; $\m{Rep}_{G}(\m{BanC}^{\m{crt}}(k))$ is the category of unitary Banach $k$-linear representations $(V,\pi)$ satisfying $\n{V} \subset \v{k}$ and contractive $k$-linear homomorphisms, and $\m{TMod}_{\m{flat,lin}}^{\m{cpt,sep}}(k^{\circ}[[G]])$ is the category of compact Hausdorff $k^{\circ}$-flat $k^{\circ}$-linear topological $k^{\circ}[[G]]$-modules and continuous $k^{\circ}[[G]]$-linear homomorphisms in Theorem \ref{int iwa - rep ban} in \S \ref{Integral-Type Duality}.

\vspace{0.1in}
The duality can be localised by tensoring $k$ to the categories. The tensor is the scalar extension of categories enriched by the monoidal categories of modules. For more detail, see the beginning of \S \ref{Duality II: Linear Case}. We show that the localisation of $\m{Rep}_{G}(\m{BanC}^{\m{crt}}(k))$ is equivalent to the category $\m{Rep}_{G}^{\m{unit}}(\m{Ban}(k))$ of unitarisable Banach $k$-linear representations and continuous $k$-linear homomorphisms in \S \ref{Banach Spaces and Bounded Homomorphisms}. We define the notions of a {\it compactly generated} locally convex $k$-vector space in Definition \ref{compactly generated}, and an {\it integralisable} compactly generated locally convex $k[[G]]$-module in Definition \ref{integralisable}. We show that the localisation of $\m{TMod}_{\m{flat,lin}}^{\m{cpt,sep}}(k^{\circ}[[G]])$ is equivalent to the category $\m{TMod}_{\m{lc,cg}}^{\m{int}}(k[[G]])$ of integralisable compactly generated locally convex $k[[G]]$-modules and continuous $k[[G]]$-linear homomorphisms in \S \ref{Compactly Generated Locally Convex Spaces}. As a consequence, we obtain duality between $\m{Rep}_{G}^{\m{unit}}(\m{Ban}(k))$ and $\m{TMod}_{\m{lc,cg}}^{\m{int}}(k[[G]])$ in Theorem \ref{iwa - rep ban}. We also verify that the irreducibility of a unitarisable Banach $k$-linear representation is equivalent to the topological simpleness of the dual $k[[G]]$-module in Corollary \ref{irr - smp}.

\vspace{0.1in}
As applications, we interpret several operations on unitarisable Banach $k$-linear representations as those on integralisable compactly generated locally convex $k[[G]]$-modules or compact Hausdorff $k^{\circ}$-flat $k^{\circ}$-linear topological $k^{\circ}[[G]]$-modules in \S \ref{Applications}. For example, we give an explicit description of the dual of the continuous induction of a unitarisable Banach $k$-linear representation of a closed subgroup in Theorem \ref{dual of the induction}. Furthermore, we calculate an explicit description of a continuous parabolic induction on $\m{GL}_n(k)$ in \S \ref{Continuous Parabolic Induction}.

\section{Convention and Definition}
\label{Convention and Definition}

Throughout this paper, we denote by $k$ a complete valuation field, by $k^{\circ} \subset k$ the ring of integral elements, and by $G$ a topological group. A complete valuation field is said to be a {\it local field} if the valuation is discrete and the residue field is a finite field. A topological group is said to be {\it profinite} if it is an inverse limit of a finite discrete group, and is said to be {\it locally profinite} if it is Hausdorff and admits an open profinite subgroup.

\subsection{Topological Module}
\label{Topological Module}

A {\it ring} is unital and associative but not necessarily commutative. For a ring $R$, we denote by $\m{Mod}(R)$ the category of left $R$-modules and $R$-linear homomorphisms. An {\it algebra} over a commutative ring is assumed to satisfy that the image of the base ring by the natural ring homomorphism lies in the centre. A {\it topological ring} is a ring $R$ endowed with a topology $\tau_R$ such that the addition
\begin{eqnarray*}
  R \times R & \to & R \\
  (r,s) & \mapsto & r + s
\end{eqnarray*}
and the multiplication
\begin{eqnarray*}
  R \times R & \to & R \\
  (r,s) & \mapsto & rs
\end{eqnarray*}
are continuous with respect to the direct product topology $\tau_R \times \tau_R$ and $\tau_R$. When we emphasise $\tau_R$, we write $R_{\tau_R}$ instead of $R$. We use a similar convention for other objects than a topological ring. A {\it topological $R_{\tau_R}$-module} is a left $R$-module $M$ endowed with a topology $\tau_M$ such that the addition
\begin{eqnarray*}
  M_{\tau_M} \times M_{\tau_M} & \to & M_{\tau_M} \\
  (m,n) & \mapsto & m + n
\end{eqnarray*}
and the scalar multiplication
\begin{eqnarray*}
  R_{\tau_R} \times M_{\tau_M} & \to & M_{\tau_M} \\
  (r,m) & \mapsto & rm
\end{eqnarray*}
are continuous. We denote by $\m{TMod}(R_{\tau_R})$ the category of topological $R_{\tau_R}$-modules and continuous $R$-linear homomorphisms. In the case when $R$ is commutative, a {\it topological $R_{\tau_R}$-algebra} is an $R$-algebra $\A$ endowed with a topology $\tau_{\A}$ such that $\A_{\tau_{\A}}$ is a topological ring and a topological $R_{\tau_R}$-module. 

\vspace{0.1in}
We regard a topological $k$-vector space as a topological $k^{\circ}$-module by restriction of scalars. The {\it adic topology} of a $k^{\circ}$-module $M$ is the topology generated by submodules of the form $aM \subset M$ for an $a \in k^{\circ} \backslash \ens{0}$. Then $M$ is a topological $k^{\circ}$-module with respect to the adic topology. If the adic topology of $M$ is Hausdorff, then we say that $M$ is adically separated.

\vspace{0.1in}
A {\it Banach $k$-vector space} is a $k$-vector space endowed with a complete norm. We denote by $\m{Ban}(k)$ (resp.\ $\m{BanC}(k)$) the category of Banach $k$-vector spaces and continuous (resp.\ contractive) $k$-linear homomorphisms. We remark that the norm of a Banach $k$-vector space is an invariant in $\m{BanC}(k)$ but not in $\m{Ban}(k)$. We denote by $\m{BanC}^{\m{crt}}(k) \subset \m{BanC}(k)$ the full subcategory of Banach $k$-vector spaces $V$ with $\n{V} \subset \v{k}$.

\begin{exm}
Let $X$ be a topological space. We denote by $\m{C}_{\m{bd}}(X,k)$ the $k$-vector space of $k$-valued bounded continuous functions on $X$, and by $\m{C}_0(X,k) \subset \m{C}_{\m{bd}}(X,k)$ the $k$-vector subspace of elements $f$ satisfying that for any $\epsilon > 0$, there is a compact subspace $Y \subset X$ such that $\v{f(x)} < \epsilon$ for any $x \in X \backslash Y$. Then  $\m{C}_{\m{bd}}(X,k)$ and $\m{C}_0(X,k)$ are Banach $k$-vector spaces with respect to the supremum norm.
\end{exm}

Note that if $X$ is a discrete space $I$, then there is a natural identification $\m{C}_{\m{bd}}(I,k) \cong (k^{\circ})^I \otimes_{k^{\circ}} k$ as $k$-vector spaces, and the induced topology on the right hand side is derived from the adic topology on the canonical lattice $(k^{\circ})^I \subset (k^{\circ})^I \otimes_{k^{\circ}} k$. Through the identification, $\m{C}_0(I,k)$ corresponds to the closure of $k^{\oplus I} \subset (k^{\circ})^I \otimes_{k^{\circ}} k$. In other words, $\m{C}_0(I,k)$ is obtained by tensoring the adic completion of $(k^{\circ})^{\oplus I}$ and $k$.

\begin{exm}
Let $M$ be a compact topological $k^{\circ}$-module. For any $v \in \m{Hom}_{\m{TMod}(k^{\circ})}(M,k)$, since $M$ is compact, the image of $M$ by $v$ is bounded. We set $\n{v} \coloneqq \sup_{m \in M} \v{v(m)}$. The map $\n{\cdot} \colon \m{Hom}_{\m{TMod}(k^{\circ})}(M,k) \to [0, \infty )$ is a complete norm of a $k$-vector space. We denote by $\m{nrm}$ the norm topology. Then $\m{Hom}_{\m{TMod}(k^{\circ})}(M,k)_{\m{nrm}}$ is a Banach $k$-vector space.
\end{exm}

\begin{exm}
Let $V_1$ and $V_2$ be Banach $k$-vector spaces. For any $A \in \m{Hom}_{\m{Ban}(k)}(V_1,V_2)$, the continuity of $A$ guarantees the boundedness of it, and hence the image of the closed unit ball by $A$ is bounded. We set
\begin{eqnarray*}
  \n{A} \coloneqq \sup \Set{\n{A(v)}}{v \in V, \n{v} \leq 1}.
\end{eqnarray*}
The map $\n{\cdot} \colon \m{Hom}_{\m{Ban}(k)}(V_1,V_2) \to [0, \infty )$ is a complete norm of a $k$-vector space. We denote by $\m{nrm}$ the norm topology. Then $\m{Hom}_{\m{Ban}(k)}(V_1,V_2)_{\m{nrm}}$ is a Banach $k$-vector space. The closed unit ball of it coincides with the $k^{\circ}$-module $\m{Hom}_{\m{BanC}(k)}(V_1,V_2)$, and the relative topology $\m{nrm}|_{\m{Hom}_{\m{BanC}(k)}(V_1,V_2)}$ coincides with the adic topology.
\end{exm}

\vspace{0.1in}
A {\it linear topological $k^{\circ}$-module} is a topological $k^{\circ}$-module admitting an open basis of $0$ consisting of open $k^{\circ}$-submodules. We denote by $\m{TMod}_{\m{flat,lin}}^{\m{cpt,sep}}(k^{\circ})$ the category of compact Hausdorff flat linear topological $k^{\circ}$-modules and continuous $k^{\circ}$-linear homomorphisms. Here the flatness of a linear topological $k^{\circ}$-module means that of the underlying $k^{\circ}$-module.

\begin{exm}
If $k^{\circ}$ is a local field, then for a set $I$, $(k^{\circ})^I$ is a compact Hausdorff flat linear topological $k^{\circ}$-module with respect to the direct product topology.
\end{exm}

Note that $(k^{\circ})^{\oplus I} \subset (k^{\circ})^I$ is a dense submodule with respect to the direct product topology. In other words, $(k^{\circ})^I$ is a compactification of $(k^{\circ})^{\oplus}$. There is duality between two topological modules $\m{C}_0(I,k)$ and $(k^{\circ})^I$ constructed from $(k^{\circ})^{\oplus I}$. See Proposition \ref{int cpt - ban} for more detail.

\vspace{0.1in}
A {\it locally convex $k$-vector space} is a topological $k$-vector space which is linear topological as a topological $k^{\circ}$-module. A {\it locally convex topology} of a $k$-vector space $V$ is a topology on $V$ for which $V$ is a locally convex $k$-vector space.

\begin{exm}
Let $M$ be a topological $k^{\circ}$-module and $W$ a Hausdorff locally convex $k$-vector space. We denote by $\m{ptw}$ (resp.\ $\m{unf}$) the topology of pointwise convergence (resp.\ uniform convergence) of $\m{Hom}_{\m{TMod}(k^{\circ})}(M,W)$. Then both of $\m{Hom}_{\m{TMod}(k^{\circ})}(M,W)_{\m{ptw}}$ and $\m{Hom}_{\m{TMod}(k^{\circ})}(M,W)_{\m{unf}}$ are Hausdorff locally convex $k$-vector spaces.
\end{exm}

\begin{exm}
Let $W$ be a $k$-vector space and $\mathscr{F}$ a family of $k$-linear homomorphisms from locally convex $k$-vector spaces to $W$. Then there is the strongest locally convex topology on $W$ for which each element of $\mathscr{F}$ is continuous.
\end{exm}

\subsection{Iwasawa Algebra}
\label{Iwasawa Algebra}

\begin{dfn}
For a topological space $X$, we denote by $\P(X)$ the set of families $P \subset 2^X$ of clopen subsets of $X$ with $\bigsqcup_{U \in P} U = X$. We endow every $P \in \P(X)$ with the discrete topology. A $P \in \P(X)$ is said to be finite if $P$ is a finite set. For a $P \in \P(X)$ and a subset $V \in X$, we set $P|_V \coloneqq \set{U \in P}{U \subset V}$. For $P_1, P_2 \in \P(X)$, we write $P_1 \leq P_2$ if for any $U \in P_1$ there exists a $V \in P_2$ such that $U \subset V$.
\end{dfn}

For $P_1, P_2 \in \P(X)$, the inequality $P_1 \leq P_2$ holds if and only if $P_1|_V \in \P(V)$ for any $V \in P_2$. The refinement $P_1 * P_2 \coloneqq \set{U \cap V}{U \in P_1, V \in P_2}$ is the least element in $\P(X)$ greater than $P_1$ and $P_2$, and hence $\P(X)$ is directed.

\begin{dfn}
For a non-empty set $I$, the set $F \subset 2^I$ of finite subsets of $I$ is directed with respect to the inclusion relation. For a topological Abelian group $M$ and a map $f \colon I \to M$, we put
\begin{eqnarray*}
  \sum_{i \in I} f(i) \coloneqq \lim_{J \in F} \sum_{i \in J} f(i).
\end{eqnarray*}
\end{dfn}

When $I = \ens{0, \ldots, n}$ for an $n \in \N$ (resp.\ $I = \N$), the value of $\sum_{i \in I}$ coincides with the usual finite sum $\sum_{i = 0}^{n}$ (resp.\ the usual infinite sum $\sum_{i = 0}^{\infty}$).

\begin{prp}
For a discrete space $X$ and an $f \in \m{C}_0(X,k^{\circ})$, the infinite sum $\sum_{x \in X} f(x)$ converges.
\end{prp}

\begin{proof}
Since $X$ is discrete, for any $\epsilon > 0$, there are at most finitely many $x \in X$ such that $\v{f(x)} > \epsilon$. In particular, $f$ is zero outside a certain countable subset $U \subset X$. Then taking a numbering on $U$, we may replace $\sum_{x \in X} f(x)$ by a sum of the form $\sum_{i = 0}^{n} f(x_i)$ or $\sum_{i = 0}^{\infty} f(x_i)$. The former case is a finite sum and hence converges. The latter sum is an infinite sum with $\lim_{i \to \infty} \v{f(x_i)} = 0$, and hence converges.
\end{proof}

Let $X$ be a topological space. For $P_1, P_2 \in \P(X)$ with $P_1 \leq P_2$, an $f \in \m{C}_0(P_1,k^{\circ})$, and a $V \in P_2$, consider the converging infinite sum $\tilde{f}(V) \coloneqq \sum_{U \in P_1|_V} f(U)$. Since $\set{P_1|_V}{V \in P_2}$ is a partition of $P_1$, the map $\tilde{f} \colon P_2 \to X \colon V \mapsto \tilde{f}(V)$ lies in $\m{C}_0(P_2,X)$. It gives a continuous $k^{\circ}$-algebra homomorphism
\begin{eqnarray*}
  \m{C}_0(P_1,k^{\circ}) & \to & \m{C}_0(P_2,k^{\circ}) \\
  f & \mapsto & \tilde{f}.
\end{eqnarray*}
Then $(\m{C}_0(P,X))_{P \in \P(X)}$ forms an inverse system.

\begin{dfn}
For a topological space $X$, set
\begin{eqnarray*}
  k^{\circ}[[X]] \coloneqq \varprojlim_{P \in \P(X)} \m{C}_0(P,X)
\end{eqnarray*}
and $k[[X]] \coloneqq k^{\circ}[[X]] \otimes_{k^{\circ}} k$. We endow $k^{\circ}[[X]]$ with the inverse limit topology $\tau_{\m{w}}$ and the weakest topology $\tau_{\m{s}}$ for which for each compact clopen subset $C \subset X$ the canonical projection $k^{\circ}[[X]] \to \m{C}_0(\ens{C,X \backslash C},k^{\circ})$ is continuous.
\end{dfn}

We do not endow $k[[G]]$ with a topology. The topology $\tau_{\m{w}}$ is stronger than or equal to $\tau_{\m{s}}$. If $X$ is compact, then $\tau_{\m{w}}$ coincides with $\tau_{\m{s}}$. Both of $\tau_{\m{w}}$ and $\tau_{\m{s}}$ are linear topologies, and hence $k^{\circ}[[X]]_{\tau_{\m{w}}}$ and $k^{\circ}[[X]]_{\tau_{\m{s}}}$ are topological $k^{\circ}$-modules. Since $\m{C}_0(P,k^{\circ})$ is Hausdorff for any $P \in \P(X)$, $\tau_{\m{w}}$ is Hausdorff. If $X$ is a disjoint union of compact clopen subsets such as a locally profinite group, then $\tau_{\m{s}}$ is also Hausdorff.

\vspace{0.1in}
The $k$-vector space $k[[X]]$ deeply related with non-Archimedean functional analysis on $X$. A {\it $k$-valued measure} on $X$ is a $k$-valued function $\mu$ on the finitely additive class consisting of clopen subsets of $X$ such that $\mu(U \sqcup V) = \mu(U) + \mu(V)$ for any disjoint clopen subsets $U,V \subset X$.

A $k$-valued measure $\mu$ on $X$ is said to be {\it bounded} (resp.\ {\it integral}) if there is an $\epsilon > 0$ such that $\v{\mu(U)} \leq \epsilon$ (resp.\ if $\v{\mu(U)} \leq 1$) for any clopen subset $U \subset X$. If $k$ is a local field, then there is a natural one-to-one correspondence between a bounded (resp.\ integral) $k$-valued measure on $X$ and a continuous (resp.\ contractive) functional on $\m{C}_{\m{bd}}(X,k)$ through integration. Note that since every bounded subset of $k$ is totally bounded, any $k$-valued bounded continuous function on $X$ is uniformly approximated by a linear combination of characteristic functions on pairwise disjoint clopen subsets.

A bounded $k$-valued measure $\mu$ on $X$ is said to be {\it normal} if the equality $\mu(V) = \sum_{U \in P} \mu(U)$ holds for any clopen subset $V \subset X$ and any $P \in \P(V)$. Then in particular $\mu(U) = 0$ for all but countably many $U \in P$ for any clopen subset $V \subset X$ and any $P \in \P(V)$. In addition if $X$ is a locally profinite group $G$, then for any $\epsilon > 0$, there is a compact clopen subset $V \subset G$ such that $\v{\mu(U)} < \epsilon$ for any clopen subset $U \subset G \backslash V$. Indeed, assume that there is an $\epsilon > 0$ such that for any compact clopen subset $V \subset X$, there is a clopen subset $U \subset G \backslash V$ with $\v{\mu(U)} \geq \epsilon$. Take a compact open subgroup $K \subset G$. By the assumption, $G$ is not compact and hence $G/K$ is an infinite set. We put $C_0 \coloneqq K$. By the assumption, there is a clopen subset $U_1 \subset G \backslash C_0$ such that $\v{\mu(U_1)} \geq \epsilon$, and by the normality of $\mu$, we may assume $U_1$ can be taken in some $C_1 \in G/K$. By the assumption again, there is a clopen subset $U_2 \subset G \backslash (C_0 \sqcup C_1)$ such that $\v{\mu(U_2)} \geq \epsilon$, and by the normality of $\mu$, we may assume $U_2$ can be taken in some $C_2 \in G/K$. Repeating this process, we obtain distinct elements $C_1, C_2, \ldots \in G/K$ and clopen subsets $U_i \subset C_i$ with $\v{\mu(U_i)} \geq \epsilon$. Since $C_1, C_2, \ldots$ are pairwise disjoint as subsets of $G$, the disjoint union $U$ of $U_1, U_2, \ldots $ is again clopen, but the infinite sum $\sum_{i \in \N} \mu(U_i)$ does not converge. This contradicts the normality of $\mu$.

For a set $I$ and a sequence $(f_i)_{i \in I}$ of $k$-valued bounded continuous functions on $X$, we say that $\sum_{i \in I} f_i$ has the normal convergence limit if for any $x \in X$, there is at most one $i \in I$ such that $f_i(x) \neq 0$. Then $\sum_{i \in I} f_i$ converges in pointwise to a $k$-valued locally bounded continuous function $f$ on $X$. The limit $f$ is bounded if and only if $(f_i)_{i \in i}$ is uniformly bounded, and then we say that $f$ is the bounded normal convergence limit. Note that any $k$-valued bounded continuous function on $X$ is uniformly approximated by a bounded locally constant function, and any bounded locally constant function is obtained as a bounded normal convergence limit of a formal linear combination of characteristic functions on pairwise disjoint clopen subsets. A normal bounded (resp.\ integral) $k$-valued measure corresponds to a continuous (resp.\ contractive) functional on $\m{C}_{\m{bd}}(X,k)$ compatible with bounded normal convergence limit through integration.

For a normal bounded (resp.\ integral) $k$-valued measure $\mu$ on $X$, the sequence $\mu_P = (\mu(U))_{U \in P}$ lies in $\m{C}_0(P,k)$ (resp.\ $\m{C}_0(P,k^{\circ})$) for any $P \in \P(X)$, and the compatible system $(\mu_P)_{P \in \P(X)}$ is an element of $k[[X]]$ (resp.\ $k^{\circ}[[X]]$). This construction is invertible, and hence we identify $k[[X]]$ (resp.\ $k^{\circ}[[X]]$) the $k$-vector space of normal bounded $k$-valued measures (resp.\ the $k^{\circ}$-module of normal integral $k$-valued measures). Each point of $X$ corresponds to a Dirac operator, and the $k^{\circ}$-submodule of $k^{\circ}[[X]]$ generated by the image of $X$ is dense with respect to $\tau_{\m{w}}$ and $\tau_{\m{s}}$. As a $k^{\circ}$-valued measure, an $x \in X$ corresponds to the truth value function on whether a clopen subset contains $x$ or not, and as a continuous functional, it corresponds to the evaluation at $x$. If $X$ is zero-dimensional, then the induced map $X \to k^{\circ}[[X]]$ is continuous. In addition if $X$ is Hausdorff, then $X \to k^{\circ}[[X]]$ is injective and the image is a linearly independent system. In particular if $X$ is a zero-dimensional Hausdorff group, then the induced homomorphism $k^{\circ}[G] \to k^{\circ}[[G]]$ is injective.

\begin{dfn}
For a commutative topological ring $R$, a bitopological $R$-algebra is an $R$-algebra $\A$ endowed with two topologies $\tau_{\m{w}}$ and $\tau_{\m{s}}$ such that $\A_{\tau_{\m{w}}}$ and $\A_{\tau_{\m{s}}}$ are topological $R$-modules, $\tau_{\m{w}}$ is stronger than $\tau_{\m{s}}$, and the multiplication
\begin{eqnarray*}
  \A_{\tau_{\m{w}}} \times \A_{\tau_{\m{w}}} & \to & \A_{\tau_{\m{s}}} \\
  (a,b) & \mapsto & ab
\end{eqnarray*}
is continuous. We call $\tau_{\m{w}}$ (resp.\ $\tau_{\m{s}}$) the weak (resp.\ strong) topology of $\A$.
\end{dfn}

When we emphasise $\tau_{\m{w}}$ and $\tau_{\m{s}}$, we write $\A_{\tau_{\m{w}},\tau_{\m{s}}}$ instead of $\A$. A topological $R$-algebra is regarded as a bitopological $R$-algebra with $\tau_{\m{w}} = \tau_{\m{s}}$ in a natural way.

\begin{dfn}
For a topological ring $R$ and bitopological $R$-algebras $\A,\A'$, a weakly continuous homomorphism from $\A$ to $\A'$ is an $R$-algebra homomorphism $\A \to \A'$ continuous with respect to their weak topologies. We denote by $\m{BTModW}(R)$ the category of bitopological $R$-algebras and weakly continuous homomorphisms.
\end{dfn}

\begin{prp}
For a locally profinite group $G$, the multiplication $k^{\circ}[G] \times k^{\circ}[G] \to k^{\circ}[G]$ extends to a unique continuous map
\begin{eqnarray*}
  k^{\circ}[[G]]_{\tau_{\m{w}}} \times k^{\circ}[[G]]_{\tau_{\m{w}}} \to k^{\circ}[[G]]_{\tau_{\m{s}}}.
\end{eqnarray*}
\end{prp}

\begin{proof}
The uniqueness, the bilinearity, the unitarity, the associativity, and the distributivity of the extended multiplication follow from the density of $k^{\circ}[G] \subset k^{\circ}[[G]]$. Therefore it suffices to show that the existence and the continuity. Let $\mu_1, \mu_2 \in k^{\circ}[[G]]$. For every compact clopen subset $O \subset G \times G$, there are an $n \in \N$ and clopen subsets $U_1^1, \ldots, U_n^1, U_1^2, \ldots, U_n^2 \subset G$ such that $\tilde{U} \cap (g_1K \times g_2K) = \bigsqcup_{i = 1}^{n} U_i^1 \times U_i^2$ because a locally profinite group is zero-dimensional. We put $(\mu_1 \times \mu_2)(O) \coloneqq \sum_{i = 1}^{n} \mu_1(U_i^1)\mu(U_i^2)$. By the finite additivity of $\mu_1$ and $\mu_2$, $(\mu_1 \times \mu_2)(O)$ depends only on $O$. In particular, $(\mu_1 \times \mu_2)(U_1 \times U_2) = \mu_1(U_1) \times \mu_2(U_2)$ for any clopen subsets $U_1, U_2 \subset X$. For a clopen subset $O \subset G \times G$ which is not necessarily compact, take a compact clopen subgroup $K \subset G$. Then $G/K \times G/K = \set{g_1K \times g_2K}{g_1, g_2 \in G}$ is a compact clopen partition of $G \times G$. We put $(\mu_1 \times \mu_2)(O) \coloneqq \sum_{C_1, C_2 \in G/K} (\mu_1 \times \mu_2)(O \cap (C_1 \times C_2))$. The infinite sum converges because $\mu_1$ and $\mu_2$ tends to $0$ at the infinity by the basic property of a normal bounded $k$-valued measure in the argument in the definition of normality. Then $(\mu_1 \times \mu_2)(O)$ is independent of the choice of $K$ by the normality of $\mu_1$ and $\mu_2$. This gives a normal integral $k$-valued measure $\mu_1 \times \mu_2$ on $G \times G$.

For each clopen subset $U \subset G$, denote by $\tilde{U} \subset G \times G$ the preimage of $U$ by the continuous multiplication $G \times G \to G$, and set $(\mu_1 * \mu_2)(U) \coloneqq (\mu_1 \times \mu_2)(\tilde{U})$. Since $\mu_1 \times \mu_2$ is a normal integral $k$-valued measure on $G \times G$, so is $\mu_1 * \mu_2$ on $G$. Thus we have constructed $\mu_1 * \mu_2 \in k^{\circ}[[G]]$. By the construction, it is easily seen that the map
\begin{eqnarray*}
  k^{\circ}[[G]]_{\tau_{\m{w}}} \times k^{\circ}[[G]]_{\tau_{\m{w}}} & \to & k^{\circ}[[G]]_{\tau_{\m{s}}} \\
  (\mu_1,\mu_2) & \mapsto & \mu_1 * \mu_2
\end{eqnarray*}
is compatible with the embedding $k^{\circ}[G] \hookrightarrow k^{\circ}[[G]]$ and the multiplication $k^{\circ}[G] \times k^{\circ}[G] \to k^{\circ}[G]$. We show the continuity. Let $\mu_1, \mu_2 \in k^{\circ}[[G]]$. For a compact clopen subset $C \subset X$ and an $\epsilon > 0$, denote by $O \subset k^{\circ}[[G]]_{\tau_{\m{s}}}$ the open neighbourhood
\begin{eqnarray*}
  \Set{\mu \in k^{\circ}[[G]]_{\tau_{\m{s}}}}{\v{\mu(C) - (\mu_1 * \mu_2)(C)}, \v{\mu(G \backslash C) - (\mu_1 * \mu_2)(G \backslash C)} < \epsilon}
\end{eqnarray*}
of $\mu_1 * \mu_2$. Note that such open neighbourhoods form an open basis. Since $G$ is locally profinite, for each $x \in C$, there is an open subgroup $K \subset G$ such that $KxK \subset C$. For any $y \in KxK$, $KxK$ contains an open subset $Ky \subset G$, and hence $KxK \subset G$ is open. Since $K$ is compact and the multiplication $G \times G \to G$ is continuous, $KxK$ is compact and hence closed in $G$ because $G$ is Hausdorff. Therefore $KxK \subset G$ is clopen. Considering all clopen subsets of $C$ of such a form, by the compactness of $C$, we obtain an $n \in \N$, $x_1, \ldots, x_n \in G$, and a compact open subgroups $K \subset G$ such that $C = \bigsqcup_{i = 1}^{n} K x_i K$. Thus the preimage of $C$ by the multiplication $G \times G \to G$ is a disjoint union of clopen subsets of $G \times G$ belonging to $K \backslash G \times G/K$, and hence the preimage of $O$ by the multiplication $k^{\circ}[[G]]_{\tau_{\m{w}}} \times k^{\circ}[[G]]_{\tau_{\m{w}}} \to k^{\circ}[[G]]_{\tau_{\m{s}}}$ contains the direct product of the open neighbourhood
\begin{eqnarray*}
  \Set{\nu \in k^{\circ}[[G]]_{\tau_{\m{w}}}}{\v{\nu(U) - \mu_1(U)} < \epsilon, {}^{\forall} U \in K \backslash G}
\end{eqnarray*}
of $\mu_1$ and the open neighbourhood
\begin{eqnarray*}
  \Set{\nu \in k^{\circ}[[G]]_{\tau_{\m{w}}}}{\v{\nu(U) - \mu_2(U)} < \epsilon, {}^{\forall} U \in G/K}
\end{eqnarray*}
of $\mu_2$. We conclude that the multiplication $k^{\circ}[[G]]_{\tau_{\m{w}}} \times k^{\circ}[[G]]_{\tau_{\m{w}}} \to k^{\circ}[[G]]_{\tau_{\m{s}}}$ is continuous.
\end{proof}

Thus $k^{\circ}[[G]]_{\tau_{\m{w}},\tau_{\m{s}}}$ is a bitopological $k^{\circ}$-algebra with respect to this multiplication. We also write $k^{\circ}[[G]]$ instead of $k^{\circ}[[G]]_{\tau_{\m{w}},\tau_{\m{s}}}$. The natural embedding $k^{\circ}[G] \hookrightarrow k^{\circ}[[G]]$ is a $k^{\circ}$-algebra homomorphism.

\begin{dfn}
For a commutative topological ring $R$ and a bitopological $R$-algebra $\A_{\tau_{\m{w}},\tau_{\m{s}}}$, a topological $\A_{\tau_{\m{w}},\tau_{\m{s}}}$-module is a left $\A$-module $M$ endowed with a topology $\tau_M$ such that $M_{\tau_M}$ is a topological $R$-module and the multiplication
\begin{eqnarray*}
  \A_{\tau_{\m{w}}} \times M_{\tau_M} & \to & M_{\tau_M} \\
  (a,m) & \mapsto & am
\end{eqnarray*}
is continuous.
\end{dfn}

\section{Duality I: Integral Case}
\label{Duality I: Integral Case}

In this section, we establish duality theory between the category of unitary Banach representations of a locally profinite group and the category of compact Hausdorff $k^{\circ}$-flat $k^{\circ}$-linear topological $k^{\circ}[[G]]$-modules. Here we deal only with contractive $k$-linear homomorphisms in the former category, and hence homsets are just $k^{\circ}$-module but not $k$-vector spaces.

\subsection{Banach Spaces and Contractive Homomorphisms}
\label{Banach Spaces and Contractive Homomorphisms}

\begin{lmm}
\label{int iwa - rep ban 1}
For a Banach $k$-vector space $V$ and a topological space $X$, the map
\begin{eqnarray*}
  (\cdot)|_G \colon \m{Hom}_{\m{TMod}(k^{\circ})} \left( k^{\circ}[[G]]_{\tau_{\m{w}}},\m{End}_{\m{BanC}(k)}(V)_{\m{ptw}} \right) & \to & \m{Hom}_{\m{Top}} \left( G,\m{End}_{\m{BanC}(k)}(V)_{\m{ptw}} \right) \\
  f & \mapsto & f|_G
\end{eqnarray*}
is bijective.
\end{lmm}

\begin{proof}
Since $k^{\circ}[G]$ is a dense $k^{\circ}$-submodule of $k^{\circ}[[G]]_{\tau_{\m{w}}}$ and $G$ generates $k^{\circ}[G]$ as a $k^{\circ}$-module, $(\cdot)|_G$ is injective. For the surjectivity, let $\varphi \colon G \to \m{End}_{\m{BanC}(k)}(V)_{\m{ptw}}$ be a continuous map. For a $v \in V$, $\varphi$ induces a continuous map
\begin{eqnarray*}
  \varphi_v \colon G & \to & V \\
  g & \mapsto & \varphi(g)(v).
\end{eqnarray*}
Since $\varphi(g)$ is contractive, $\varphi(g)(v)$ is contained in the closed disc $B_v \subset V$ of radius $\n{v}$ centred at $0$ for any $g \in G$. Consider the $k^{\circ}$-linear extension $f_v \colon k^{\circ}[G] \to V$ of $\varphi_v$. Since $B_v \subset V$ is a $k^{\circ}$-submodule, $f_v(k^{\circ}[G])$ is contained in $B_v$. We verify the continuity of $f_v$ with respect to the relative topology of $k^{\circ}[G] \subset k^{\circ}[[G]]_{\tau_{\m{w}}}$. Let $W \subset V$ be an open subgroup. Denote by $P$ the partition of $G$ given as the preimage of the partition $V/W = \set{v + W}{v \in V}$ of $V$. Take a $u \in k^{\circ} \backslash \ens{0}$ such that $uv \in W$. Then $u \varphi(g)(v) \in W$ because $\varphi(g)$ is contractive for any $g$. Therefore $f_v^{-1}(W)$ contains the open ideal of $k^{\circ}[G]$ which is the preimage of $u \m{C}_0(P,k^{\circ})$ by the composite $k^{\circ}[G] \hookrightarrow k^{\circ}[[G]]_{\tau_{\m{w}}} \twoheadrightarrow \m{C}_0(P,k^{\circ})$. This gives the continuity of $f_v$. Note that any continuous group homomorphism between topological groups is uniformly continuous with respect to the canonical uniform structures, and hence so is $f_v$. Therefore $f_v$ continuously extends to $k^{\circ}[[G]]_{\tau_{\m{w}}}$ by the completeness of $V$. We also denote the extension by the same symbol $f_v$. By the density of $k^{\circ}[G] \subset k^{\circ}[[G]]_{\tau_{\m{w}}}$, the extension is also $k^{\circ}$-linear. Since $B_v \subset V$ is closed, $f_v(k^{\circ}[[G]])$ is contained in $B_v$. For a $\mu \in k^{\circ}[[G]]$, consider a map
\begin{eqnarray*}
  f(\mu) \colon V & \to & V \\
  v & \mapsto & f_v(\mu).
\end{eqnarray*}
Since the map $V \to \m{Hom}_{\m{Top}}(G,V) \colon v \mapsto \varphi_v$ is $k^{\circ}$-linear, $f(\mu)$ is $k$-linear. Moreover, $f(\mu)$ is contractive because $f_v(\mu)$ is contained in $B_v$ for any $v \in V$. This gives a map
\begin{eqnarray*}
  f \colon k^{\circ}[[G]]_{\tau_{\m{w}}} & \to & \m{End}_{\m{BanC}(k)}(V)_{\m{ptw}} \\
  \mu & \mapsto & f(\mu).
\end{eqnarray*}
Since $f_v$ is $k^{\circ}$-linear for each $v \in V$, so is $f$. The continuity of $f_v$ for each $v \in V$ guarantees that of $f$ with respect to the topology of pointwise convergence. Thus $f \in \m{Hom}_{\m{TMod}(k^{\circ})}(k^{\circ}[[G]]_{\tau_{\m{w}}},\m{End}_{\m{BanC}(k)}(V)_{\m{ptw}})$, and is an extension of $\varphi$ by the construction. We conclude that $(\cdot)|_G$ is surjective.
\end{proof}

\begin{prp}
\label{int iwa - rep ban 2}
Suppose that $G$ is locally profinite. For a Banach $k$-vector space $V$, the map
\begin{eqnarray*}
  (\cdot)|_G \colon \m{Hom}_{\m{BTAlgW}(k^{\circ})} \left( k^{\circ}[[G]],\m{End}_{\m{BanC}(k)}(V)_{\m{ptw}} \right) & \to & \m{Hom}_{\m{TGrp}} \left( G,\m{Aut}_{\m{BanC}(k)}(V)_{\m{ptw}} \right) \\
  f & \mapsto & f|_G
\end{eqnarray*}
is bijective.
\end{prp}

\begin{proof}
The map is well-defined because $G$ is contained in $k^{\circ}[[G]]^{\times}$. The injectivity follows from Lemma \ref{int iwa - rep ban 1}. For a continuous group homomorphism $\varphi \colon G \to \m{Aut}_{\m{BanC}(k)}(V)_{\m{ptw}}$, take a unique extension $f \colon k^{\circ}[[G]]_{\tau_{\m{w}}} \to \m{End}_{\m{BanC}(k)}(V)_{\m{ptw}}$ as a continuous $k^{\circ}$-linear homomorphism. By the universality of the group algebra $k^{\circ}[G]$, the restriction $f|_{k^{\circ}[G]}$ is a $k^{\circ}$-algebra homomorphism, and hence so is $f$ by the density argument. As a consequence, $(\cdot)|_G$ is surjective.
\end{proof}

\begin{lmm}
\label{rep - ptw}
For a Banach $k$-vector space $V$ and a topological space $X$, a map $\tilde{\rho} \colon X \to \m{End}_{\m{BanC}(k)}(V)_{\m{ptw}}$ is continuous if and only if the induced action
\begin{eqnarray*}
  \rho \colon X \times V & \to & V \\
  (x,v) & \mapsto & \tilde{\rho}(x)(v)
\end{eqnarray*}
is continuous.
\end{lmm}

\begin{proof}
Suppose that $\rho$ is continuous. For every $v \in V$, the restriction
\begin{eqnarray*}
  \rho|_{X \times \ens{v}} \colon X & \to & V \\
  x & \mapsto & \rho(x,v)
\end{eqnarray*}
is continuous. This implies the continuity of $\tilde{\rho}$ by the definition of the topology of pointwise convergence.

Suppose that $\tilde{\rho}$ is continuous. Let $x \in X$ and $v \in V$. By the continuity of $\tilde{\rho}$, for any open ball $B \subset V$ containing $\tilde{\rho}(x)(v)$, there exists an open neighbourhood $U \subset X$ of $x$ such that $\tilde{\rho}(u)(v) \in B$ for any $u \in U$. Since $\tilde{\rho}(x)$ is contractive, there exists an open ball $B' \subset V$ containing $v$ such that $\tilde{\rho}(u)(v') \in B$ for any $u \in U$ and $v' \in B'$. Therefore $U \times B'$ is an open neighbourhood of $(x,v) \in X \times V$ contained in $\rho^{-1}(B)$. It implies that $\rho$ is continuous.
\end{proof}

Thus applying Lemma \ref{rep - ptw} to the case when $X=G$, giving a strongly continuous unitary action of $G$ on $V$ is equivalent to giving a continuous group homomorphism $\tilde{\rho} \colon G \to \m{Aut}_{\m{BanC}(k)}(V)_{\m{ptw}}$.

\begin{dfn}
A unitary Banach $k$-linear representation of $G$ is a pair $(V,\pi)$ of a Banach $k$-vector space $V$ and a continuous map $\pi \colon G \times V \to V$ such that for each $g$, the map $\pi(g,\cdot) \colon V \to V \colon v \mapsto \pi(g,v)$ is a contractive $k$-linear automorphism and the induced map $\tilde{\pi} \colon G \to \m{Aut}_{\m{BanC}(k)}(V)$ is a group homomorphism. A unitary Banach $k$-linear representation $(V,\pi)$ of $G$ is said to be strictly Cartesian if $\n{V} \subset \v{k}$. Denote by $\m{Rep}_{G}(\m{BanC}(k))$ the category of unitary Banach $k$-linear representations of $G$ and contractive $G$-equivariant $k$-linear homomorphisms, and by $\m{Rep}_{G}(\m{BanC}^{\m{crt}}(k)) \subset \m{Rep}_{G}(\m{BanC}(k))$ the full subcategory of strictly Cartesian unitary Banach $k$-linear representations $(V,\pi)$ of $G$.
\end{dfn}

\begin{dfn}
An integral Banach $k[[G]]$-module is a $k[[G]]$-module $V$ endowed with a complete norm of a $k$-vector space such that the closed unit ball of $V$ is a topological $k^{\circ}[[G]]$-module with respect to the norm topology. An integral Banach $k[[G]]$-module is said to be strictly Cartesian if the underlying Banach $k$-vector space $V$ satisfies $\n{V} \subset \v{k}$. Denote by $\m{BanC}(k^{\circ}[[G]])$ the category of integral Banach $k[[G]]$-modules and contractive $k[[G]]$-linear homomorphisms, and by $\m{BanC}^{\m{crt}}(k^{\circ}[[G]]) \subset \m{BanC}(k^{\circ}[[G]])$ the full subcategory of strictly Cartesian integral Banach $k[[G]]$-module.
\end{dfn}

For a topological $k^{\circ}[[G]]$-module $M$, we denote by$M|_G$ the $k^{\circ}$-linear representation $(M,\pi|_{G \times M})$ of $G$, where $\pi$ is the scalar multiplication $k^{\circ}[[G]] \times M \to M$.

\begin{prp}
\label{int iwa - rep ban 3}
Suppose that $G$ is locally profinite. The covariant functor
\begin{eqnarray*}
  (\cdot)|_G \colon \m{BanC}(k^{\circ}[[G]]) & \to & \m{Rep}_{G}(\m{BanC}(k)) \\
  V & \rightsquigarrow & V|_G
\end{eqnarray*}
is an equivalence of categories.
\end{prp}

\begin{proof}
The fullness and the faithfulness is obvious by the density of $k^{\circ}[G] \subset k^{\circ}[[G]]_{\tau_{\m{w}}}$. Let $(V,\varpi)$ be an object of $\m{Rep}_{G}(\m{BanC}(k))$. The action $\varpi \colon G \times V \to V$ induces a continuous group homomorphism $\tilde{\varpi} \colon G \to \m{Aut}_{\m{BanC}(k)}(V)_{\m{ptw}}$ by Lemma \ref{rep - ptw}, and $\tilde{\varpi}$ extends to a weakly continuous $k^{\circ}$-algebra homomorphism $\tilde{\pi} \colon k^{\circ}[[G]] \to \m{End}_{\m{BanC}(k)}(V)_{\m{ptw}}$ by Proposition \ref{int iwa - rep ban 2}. Moreover, $\tilde{\pi}$ induces an continuous contractive action $\pi \colon k^{\circ}[[G]]_{\tau_{\m{w}}} \times V \to V$ again by Lemma \ref{rep - ptw}. Then $(V,\pi)$ is a topological $k^{\circ}[[G]]$-module, and the equality $\pi|_{G \times V} = \varpi$ holds by the construction. The construction of $(V,\pi)$ is natural on $(V,\varpi)$, and hence by the full faithfulness of $(\cdot)|_G$, it gives a contravariant functor
\begin{eqnarray*}
  \m{Rep}_{G}(\m{BanC}(k)) \to \m{BanC}(k^{\circ}[[G]])
\end{eqnarray*}
which is the inverse of $(\cdot)|_G$.
\end{proof}

Restricting the equivalence to the full subcategories, we obtain the following equivalence.

\begin{crl}
\label{int iwa - rep ban 4}
Suppose that $G$ is locally profinite. The covariant functor
\begin{eqnarray*}
  (\cdot)|_G \colon \m{BanC}^{\m{crt}}(k^{\circ}[[G]]) & \to & \m{Rep}_{G} \left( \m{BanC}^{\m{crt}}(k) \right) \\
  V & \rightsquigarrow & V|_G
\end{eqnarray*}
is an equivalence of categories.
\end{crl}

\subsection{Compact Linear Topological Modules}
\label{Compact Linear Topological Modules}

\begin{lmm}
\label{compact - adic}
For a compact linear topological $k^{\circ}$-module $M$ and an open neighbourhood $N \subset M$ of $0$, there is an $a \in k^{\circ} \backslash \ens{0}$ such that $aM \subset N$.
\end{lmm}

\begin{proof}
We may assume $N \subset M$ is an open $k^{\circ}$-submodule. Since $M$ is compact, $M/N$ is a finite $k^{\circ}$-module, and hence there is an $a \in k^{\circ}$ such that $a(M/N) = O$. It implies $aM \subset N$.
\end{proof}

Thus any compact linear topology on a $k^{\circ}$-module is weaker than the adic topology. In particular, a $k^{\circ}$-module which admits a compatible compact Hausdorff linear topology is adically separated.

\begin{lmm}
\label{int iwa - rep cpt 1}
For a compact Hausdorff linear topological $k^{\circ}$-module $M$, the map
\begin{eqnarray*}
  (\cdot) |_G \colon \m{Hom}_{\m{TMod}(k^{\circ})} \left( k^{\circ}[[G]],\m{End}_{\m{TMod}(k^{\circ})}(M)_{\m{unf}} \right) & \to & \m{Hom}_{\m{Top}} \left( G,\m{End}_{\m{TMod}(k^{\circ})}(M)_{\m{unf}} \right) \\
  f & \mapsto & f|_G
\end{eqnarray*}
is bijective.
\end{lmm}

\begin{proof}
The injectivity is verified in the same way as in Lemma \ref{int iwa - rep ban 1}. For the surjectivity, let $\varphi \colon G \to \m{End}_{\m{TMod}(k^{\circ})}(M)_{\m{unf}}$ be a continuous map. Denote by $f \colon k^{\circ}[G] \to \m{End}_{\m{TMod}(k^{\circ})}(M)_{\m{unf}}$ the $k^{\circ}$-linear extension of $\varphi$. Let $(\mu_i)_{i \in I}$ be a Cauchy net in $k^{\circ}[G]$ with respect to the relative uniform structure of $k^{\circ}[G] \subset k^{\circ}[[G]]_{\tau_{\m{w}}}$. By the continuity of $\varphi$, for an open $k^{\circ}$-module $N \subset M$, the partition $P$ of $G$ given as the preimage of the partition $M/N = \set{m + N}{m \in M}$ of $M$ consists of pairwise disjoint clopen subsets. By Lemma \ref{compact - adic}, there is an $a \in k^{\circ} \backslash \ens{0}$ such that $aM \subset N$. Denote by $O \subset k^{\circ}[G]$ the preimage of the open ideal $a \m{C}_0(P,k^{\circ})$ by the composite $k^{\circ}[G] \hookrightarrow k^{\circ}[[G]]_{\tau_{\m{w}}} \twoheadrightarrow \m{C}_0(P,k^{\circ})$. Since $(\mu_i)_{i \in I}$ is a Cauchy net, there is an $i_0 \in I$ such that $\mu_i - \mu_{i'} \in O$ for any $i,i' \geq i_0$. Then $f(\mu_i)(m) - f(\mu_{i'})(m) = f(\mu_i - \mu_{i'})(m) \in N$ for any $i,i' \geq i_0$ and $m \in M$, and hence $(f(\mu_i))_{i \in I}$ is a Cauchy net in $\m{End}_{\m{TMod}(k^{\circ})}(M)_{\m{unf}}$. Since $\m{End}_{\m{TMod}(k^{\circ})}(M)_{\m{unf}}$ is complete, $(f(\mu_i))_{i \in I}$ admits a unique limit. The density of $k^{\circ}[G] \subset k^{\circ}[[G]]_{\tau_{\m{w}}}$ guarantees that $f$ uniquely extends to a continuous $k^{\circ}$-linear homomorphism $f \colon k^{\circ}[[G]]_{\tau_{\m{w}}} \to \m{End}_{\m{TMod}(k^{\circ})}(M)_{\m{unf}}$ satisfying $f|_G = \varphi$.
\end{proof}

\begin{prp}
\label{int iwa - rep cpt 2}
Suppose that $G$ is locally profinite. For a compact Hausdorff linear topological $k^{\circ}$-module $M$, the map
\begin{eqnarray*}
  (\cdot) |_G \colon \m{Hom}_{\m{BTAlgW}(k^{\circ})} \left( k^{\circ}[[G]],\m{End}_{\m{TMod}(k^{\circ})}(M)_{\m{unf}} \right) & \to & \m{Hom}_{\m{TGrp}} \left( G,\m{Aut}_{\m{TMod}(k^{\circ})}(M)_{\m{unf}} \right) \\
  f & \mapsto & f|_G
\end{eqnarray*}
is bijective.
\end{prp}

\begin{proof}
The assertion is verified by Lemma \ref{int iwa - rep cpt 1} in the same way as in Proposition \ref{int iwa - rep ban 2}.
\end{proof}

\begin{lmm}
\label{rep - unf}
For a topological space $X$ and a compact uniform space $Y$, a map $\rho \colon X \to \m{End}_{\m{Unf}}(Y)_{\m{unf}}$ is continuous if and only if the induced action
\begin{eqnarray*}
  \tilde{\rho} \colon X \times Y & \to & Y \\
  (x,y) & \mapsto & \rho(x)(y)
\end{eqnarray*}
is continuous.
\end{lmm}

\begin{proof}
Suppose that $\rho$ is continuous. Let $x \in X$ and $y \in Y$. For an open neighbourhood $O \subset Y$ of $\tilde{\rho}(x,y)$, take an entourage $E \subset Y \times Y$ such that $E_{\tilde{\rho}(x,y)} = \set{y' \in Y}{(\tilde{\rho}(x,y),y') \in E}$ is contained in $O$. Then there is an entourage $E' \subset Y \times Y$ such that the condition $(y',y'''), (y''',y'') \in E'$ implies $(y',y'') \in E$. By the uniform continuity of $\rho(x)$, there exists an entourage $E'' \subset Y \times Y$ such that the condition $(y',y'') \in E''$ implies $(\rho(x)(y'),\rho(x)(y'')) \in E'$. By the continuity of $\rho$, there exists an open neighbourhood $U \subset X$ of $x$ such that $(\rho(u)(y'),\rho(x)(y')) \in E'$ for any $u \in U$ and $y' \in Y$. Then for any $u \in U$ and any $y' \in E''_y = \set{y' \in Y}{(y,y') \in E''}$, we have $(\tilde{\rho}(x,y),\tilde{\rho}(x,y')), (\tilde{\rho}((x,y'),\tilde{\rho}(u,y')) \in E'$, and hence $(\tilde{\rho}(x,y),\tilde{\rho}(u,y')) \in E$. This shows that $U \times E''_y \subset X \times Y$ is an open neighbourhood of $(x,y)$ contained in $O$. As a consequence, $\tilde{\rho}$ is continuous.

Suppose that $\tilde{\rho}$ is continuous. Let $x \in X$. For each entourage $E \subset Y \times Y$, set $\Fil_E \coloneqq \set{A \in \m{End}_{\m{Unf}}(Y)_{\m{unf}}}{(\rho(x)(y),A(y)) \in E, {}^{\forall} y \in Y}$. Note that the collection of subsets of the form $\Fil_E$ forms an open basis of $\rho(x)$. For an entourage $E \subset Y \times Y$, take an entourage $E' \subset Y \times Y$ such that the condition $(y,y'), (y,y'') \in E'$ implies $(y',y'') \in E$. By the continuity of $\tilde{\rho}$, for each $y \in Y$, there exist open neighbourhoods $U \subset X$ and $O \subset Y$ of $x$ and $y$ such that $(\tilde{\rho}(x,y),\tilde{\rho}(u,y')) \in E'$ for any $u \in U$ and any $y' \in O$. Consider the set $\Fil$ of such a $(y,U,O)$ for every $y \in Y$. Since $Y$ is compact, there exist an $n \in \N$ and $(y_1,U_1,O_1), \ldots, (y_n,U_n,O_n) \in \Fil$ such that $Y = O_1 \cup \cdots \cup O_n$. Set $U \coloneqq U_1 \cap \cdots \cap U_n$. Then for any $u \in U$ and any $y' \in Y$, taking an $i = 1, \ldots, n$ with $y' \in O_i$, we have $(\tilde{\rho}(x,y_i),\tilde{\rho}(u,y')) \in E'$. Since $x \in U$, it implies $(\tilde{\rho}(x,y_i), \tilde{\rho}(x,y')) \in E'$ and hence $(\tilde{\rho}(x,y'), \tilde{\rho}(u,y')) \in E$ for any $y' \in Y$. This shows that $U \subset X$ is an open neighbourhood of $x$ contained in $\rho^{-1}(\Fil_E)$. As a consequence, $\rho$ is continuous.
\end{proof}

Thus applying these results to the case when $X=G$ and $Y$ is a compact Abelian group $M$, giving a continuous additive action of $G$ on $M$ is equivalent to giving a continuous group homomorphism $\rho \colon G \to \m{Aut}_{\m{TGrp}}(M)_{\m{unif}}$. Moreover, if $M$ is a compact $k^{\circ}$-module, giving a continuous $k^{\circ}$-linear action of $G$ on $M$ is equivalent to giving a continuous homomorphism $\rho \colon G \to \m{Aut}_{\m{TMod}(k^{\circ})}(M)_{\m{unf}}$.

\begin{dfn}
A topological $k^{\circ}[[G]]$-module is said to be $k^{\circ}$-linear topological if its underlying topological $k^{\circ}$-module is linear topological. We denote by $\m{TMod}_{\m{flat,lin}}^{\m{cpt,sep}}(k^{\circ}[[G]])$ the category of compact Hausdorff $k^{\circ}$-flat $k^{\circ}$-linear topological $k^{\circ}[[G]]$-modules and continuous $k^{\circ}[[G]]$-linear homomorphism.
\end{dfn}

\begin{prp}
\label{int iwa - rep cpt 3}
Suppose that $G$ is locally profinite. The covariant functor
\begin{eqnarray*}
  (\cdot) |_G \colon \m{TMod}_{\m{flat,lin}}^{\m{cpt,sep}}(k^{\circ}[[G]]) & \to & \m{Rep}_{G} \left( \m{TMod}_{\m{flat,lin}}^{\m{cpt,sep}}(k^{\circ}) \right) \\
  M & \rightsquigarrow & M|_G
\end{eqnarray*}
is an equivalence of categories.
\end{prp}

\begin{proof}
The assertion is verified by Proposition \ref{int iwa - rep cpt 2} and Lemma \ref{rep - unf} in the same way as Proposition \ref{int iwa - rep ban 3}. Note that the equivalence holds without restricting to flat $k^{\circ}$-modules, and the equivalence respects the flatness.
\end{proof}

\subsection{Integral-Type Duality}
\label{Integral-Type Duality}

In this section, suppose that $k$ is a local field and $G$ is locally profinite. Then $k^{\circ}$ is compact. For a Banach $k$-vector space $V$, we set $V^D \coloneqq \m{Hom}_{\m{BanC}(k)}(V,k)_{\m{ptw}}$. Since $k^{\circ}$ is compact, $V^D$ is a compact Hausdorff flat linear topological $k^{\circ}$-module. For Banach $k$-vector spaces $V_1$ and $V_2$, every morphism $A \colon V_1 \to V_2$ in $\m{BanC}(k)$ induces a $k^{\circ}$-linear homomorphism $A^D \colon V_2^D \to V_1^D \colon m \mapsto m \circ A$. Then $A^D$ is continuous. Thus the correspondences $V \rightsquigarrow V^D$ and $A \rightsquigarrow A^D$ give a contravariant functor $D \colon \m{BanC}(k) \to \m{TMod}_{\m{flat,lin}}^{\m{cpt,sep}}(k^{\circ})$. Restricting it, we obtain a contravariant functor $D \colon \m{BanC}^{\m{crt}}(k) \to \m{TMod}_{\m{flat,lin}}^{\m{cpt,sep}}(k^{\circ})$. On the other hand, for a compact topological $k^{\circ}$-module $M$, we set $M^D \coloneqq \m{Hom}_{\m{TMod}(k^{\circ})}(M,k)_{\m{nrm}}$. For compact topological $k^{\circ}$-modules $M_1$ and $M_2$, every morphism $B \colon M_1 \to M_2$ in $\m{TMod}(k^{\circ})$ induces a $k$-linear homomorphism $B^D \colon M_2^D \to M_1^D \colon v \mapsto v \circ B$. Then $B^D$ is continuous. Thus the correspondences $M \rightsquigarrow M^D$ and $B \rightsquigarrow B^D$ give a contravariant functor $D \colon \m{TMod}_{\m{flat,lin}}^{\m{cpt,sep}}(k^{\circ}) \to \m{BanC}(k)$. If the valuation of $k$ is trivial, discrete, or surjective, then the essential image of $D$ is contained in $\m{BanC}^{\m{crt}}(k)$, and hence $D$ induces a contravariant functor $D \colon \m{TMod}_{\m{flat,lin}}^{\m{cpt,sep}}(k^{\circ}) \to \m{BanC}^{\m{crt}}(k)$.

\begin{prp}
\label{int cpt - ban}
The contravariant functors
\begin{eqnarray*}
  D \colon \m{BanC}^{\m{crt}}(k) \to \m{TMod}_{\m{flat,lin}}^{\m{cpt,sep}}(k^{\circ})
\end{eqnarray*}
and
\begin{eqnarray*}
  D \colon \m{TMod}_{\m{flat,lin}}^{\m{cpt,sep}}(k^{\circ}) \to \m{BanC}^{\m{crt}}(k)
\end{eqnarray*}
are the inverses of each other.
\end{prp}

Note that the assertion is verified in the proof of \cite{ST} Theorem 1.2 when $k$ is of characteristic $0$. It is obvious that the same proof is valid even if the assumption of the characteristic is removed, but we give an alternative proof for convenience.

\begin{proof}
By \cite{Sch} Remark 10.2, any Banach $k$-vector space $V$ with $\n{V} \subset \v{k}$ is isometrically isomorphic to $\m{C}_0(I,k)$ for some discrete space $I$. By \cite{SGA3-1} Expose VII${}_B$ 0.3.8.\ Corollaire, any compact Hausdorff flat linear topological $k^{\circ}$-module is isomorphic to $(k^{\circ})^I$ in $\m{TMod}(k^{\circ})$ for some set $I$. It is easy to see that the canonical pairing
\begin{eqnarray*}
  \m{C}_0(I,k) \otimes_{k^{\circ}} (k^{\circ})^I & \to & k \\
  (f_i)_{i \in I} \otimes (a_i)_{i \in I} & \mapsto & \sum_{i \in I} a_i f_i
\end{eqnarray*}
gives natural isomorphisms $\m{C}_0(I,k)^D \cong (k^{\circ})^I$ in $\m{TMod}_{\m{flat,lin}}^{\m{cpt,sep}}(k^{\circ})$ and $((k^{\circ})^I)^D \cong \m{C}_0(I,k)$ in $\m{BanC}^{\m{crt}}(k)$.

For a Banach $k$-vector space $V$ with $\n{V} \subset \v{k}$ and a discrete space $I$ with an isometric isomorphism $U \colon V \stackrel{\sim}{\to} \m{C}_0(I,k)$ in $\m{BanC}^{\m{crt}}(k)$, the induced isomorphism
\begin{eqnarray*}
\begin{CD}
  V @>{U}>> \m{C}_0(I,k) @>{\sim}>> \left( (k^{\circ})^I \right)^D @>{(U^{-1})^{DD}}>> V^{DD}
\end{CD}
\end{eqnarray*}
is independent of $I$ and $U$, and gives a natural equivalence $\m{id}_{\m{BanC}^{\m{crt}}(k)} \stackrel{\sim}{\to} DD$. A natural equivalence $\m{id}_{\m{TMod}_{\m{flat,lin}}^{\m{cpt,sep}}(k^{\circ})} \stackrel{\sim}{\to} DD$ is also constructed in a similar way.
\end{proof}

The duality functor is easily extended to the categories of representations. Let $(M,\rho)$ be an object of $\m{Rep}_{G}( \m{TMod}_{\m{flat,lin}}^{\m{cpt,sep}}(k^{\circ}))$. For a $(g,v) \in G \times M^D \times M$, we define a map $\rho^D(g,v) \colon M \to k$ by setting $\rho^D(g,v)(m) \coloneqq v(\rho(g^{-1},m))$ for each $m \in M$. Then the map$\rho^D(g,v) \colon M \to k \colon m \mapsto \rho^D(g,v)(m)$ is a continuous $k^{\circ}$-linear homomorphism, and 
\begin{eqnarray*}
  \rho^D \colon G \times M^D & \to & M^D \\
  (g,v) & \mapsto & \rho^D(g,v)
\end{eqnarray*}
is a $k$-linear action of $G$ such that each $g \in g$ corresponds to an isometric automorphism of $M^D$. Since $\m{End}_{\m{BanC}(k)}(M^D)$ is equicontinuous by Banach--Steinhaus theorem (\cite{Sch} Corollary 6.16), the separated continuity of $\rho^D$ guarantees the continuity of it.

\begin{lmm}
\label{int rep cpt - rep ban}
The contravariant functor
\begin{eqnarray*}
  D \colon \m{Rep}_{G} \left( \m{TMod}_{\m{flat,lin}}^{\m{cpt,sep}}(k^{\circ}) \right) & \to & \m{Rep}_{G} \left( \m{BanC}^{\m{crt}}(k) \right) \\
  (M, \rho) & \rightsquigarrow & (M^D, \rho^D)
\end{eqnarray*}
is an equivalence of categories.
\end{lmm}

\begin{proof}
The forgetful functors 
\begin{eqnarray*}
  \m{Rep}_{G} \left( \m{TMod}_{\m{flat,lin}}^{\m{cpt,sep}}(k^{\circ}) \right) & \to & \m{TMod}_{\m{flat,lin}}^{\m{cpt,sep}}(k^{\circ}) \\
  (M, \rho) & \rightsquigarrow & M
\end{eqnarray*}
and
\begin{eqnarray*}
  \m{Rep}_{G} \left( \m{BanC}^{\m{crt}}(k) \right) & \to & \m{BanC}^{\m{crt}}(k) \\
  (V, \pi) & \rightsquigarrow & V
\end{eqnarray*}
are faithful, and $D$ is compatible with the duality functor
\begin{eqnarray*}
  \m{TMod}_{\m{flat,lin}}^{\m{cpt,sep}}(k^{\circ}) & \to & \m{BanC}^{\m{crt}}(k) \\
  M & \rightsquigarrow & M^D.
\end{eqnarray*}
We distinguish the last functor from $D$ by the symbol $D'$. Then $D'$ is an equivalence by the proof of Theorem 1.2 in \cite{ST}, and hence $D$ is faithful. For the fullness, let $(M_1, \rho_1)$ and $(M_2,\rho_2)$ be an object of $\m{Rep}_{G}(\m{TMod}_{\m{flat,lin}}^{\m{cpt,sep}}(k^{\circ}))$, and $A \colon (M_2^D,\rho_2^D) \to (M_1^D,\rho_1^D)$ a morphism in $\m{Rep}_{G}(\m{BanC}^{\m{crt}}(k))$. The underlying contractive homomorphism $A \colon M_2^D \to M_1^D$ corresponds to a unique continuous $k^{\circ}$-linear homomorphism $\alpha \colon M_1 \to M_2$ because of the fully faithfulness of $D'$. Then $\alpha$ coincides with the composite
\begin{eqnarray*}
  M_1 \cong M_1^{DD} & \to & M_2^{DD} \cong M_2 \\
  h & \mapsto & h \circ A,
\end{eqnarray*}
and hence is $G$-equivariant. Thus $D$ is full.

Let $(V,\pi)$ be an object of $\m{Rep}_{G}(\m{BanC}^{\m{crt}}(k))$. The unitary action $\pi \colon G \times V \to V$ gives $\tilde{\varpi} \colon k^{\circ}[[G]] \to \m{End}_{\m{BanC}^{\m{crt}}(k)}(V)_{\m{ptw}}$ by Proposition \ref{int iwa - rep ban 2} and Lemma \ref{rep - ptw}. Set $M \coloneqq V^D$. By the proof of Proposition 1.6 in \cite{ST}, the duality map
\begin{eqnarray*}
  \m{End}_{\m{TMod}(k^{\circ})}(M)_{\m{unf}} \to \m{End}_{\m{BanC}^{\m{crt}}(k)}(V)_{\m{ptw}}
\end{eqnarray*}
is a homeomorphic anti-multiplicative isomorphism. Note that $k$ is assumed to be of characteristic $0$ in Proposition 1.6 in \cite{ST} but the proof is valid for the case when $k$ is of positive characteristic. Composing  $\tilde{\varphi}$ with the homeomorphic anti-multiplicative automorphism on $k^{\circ}[[G]]$ associated to the inverse $G \to G^{\m{op}} \colon g \mapsto g^{-1}$, we obtain a continuous homomorphism $\tilde{\rho} \colon k^{\circ}[[G]] \to \m{End}_{\m{TMod}(k^{\circ})}(M)$. By Proposition \ref{int iwa - rep cpt 2} and Lemma \ref{rep - unf}, $\tilde{\rho}$ corresponds to $\rho \colon G \times M \to M$. Then by the construction, there is a natural isomorphism $(M, \rho)^D \cong (V, \pi)$. The construction of $(M, \rho)$ is natural on $(V, \pi)$, and hence by the full faithfulness of $D$, it gives a contravariant functor
\begin{eqnarray*}
  \m{Rep}_{G} \left( \m{BanC}^{\m{crt}}(k) \right) \to \m{Rep}_{G} \left( \m{TMod}_{\m{flat,lin}}^{\m{cpt,sep}}(k^{\circ}) \right),
\end{eqnarray*}
which is the inverse of $D$.
\end{proof}

\begin{thm}
\label{int iwa - rep ban}
Suppose that $k$ is a local field and $G$ is locally profinite. There is a contravariant equivalence
\begin{eqnarray*}
  D \colon \m{TMod}_{\m{flat,lin}}^{\m{cpt,sep}}(k^{\circ}[[G]]) \to \m{Rep}_{G} \left( \m{BanC}^{\m{crt}}(k) \right)
\end{eqnarray*}
of $\m{Mod}(k^{\circ})$-enriched categories.
\end{thm}

For the notion of $\m{Mod}(k^{\circ})$-enriched category, see the beginning of \S \ref{Duality II: Linear Case}.

\begin{proof}
The equivalence is given as the composite of the equivalences in Proposition \ref{int iwa - rep cpt 3} and Lemma \ref{int rep cpt - rep ban}, and hence is an equivalence.
\end{proof}

\section{Duality II: Linear Case}
\label{Duality II: Linear Case}

In this section, we localise the results in \S \ref{Duality I: Integral Case} by tensoring $k$ to the categories in the following sense. For a ring $R$, we denote by $\m{Mod}(R)$ the category of left $R$-modules and $R$-linear homomorphisms. If $R$ is commutative, then $\m{Mod}(R)$ is a monoidal category with respect to the tensor product bifunctor $\otimes_R$ and the natural equivalences $R \otimes_R (\cdot) \cong \m{id}_{\m{Mod}(R)} \cong (\cdot) \otimes_R R$. We call a $\m{Mod}(R)$-enriched category an {\it $R$-category} for short. An {\it $R$-linear functor} is an enriched functor between $R$-categories, i.e.\ a functor between the underlying categories such that the associated maps between homsets are $R$-linear homomorphisms. An {\it $R$-linear natural transform} is a natural transform between $R$-linear functors. We denote by $\mathbb{M} \m{od}(R)$ the $2$-category of $R$-categories, $R$-linear functors, and $R$-linear natural transforms. For a commutative $R$-algebra $R'$, $\m{Mod}(R')$ is $\m{Mod}(R)$-enriched with respect to the natural monoidal functor associated to the structure morphism $R \to R'$. Therefore every $R'$-category is naturally $\m{Mod}(R)$-enriched, and every $R'$-linear functor is $R$-linear. The monoidal functor
\begin{eqnarray*}
  (\cdot) \otimes_R R' \colon \m{Mod}(R) & \to & \m{Mod}(R') \\
  M & \rightsquigarrow & M \otimes_R R'
\end{eqnarray*}
is $R$-linear, and hence it induces a $2$-functor
\begin{eqnarray*}
  (\cdot) \otimes_R R' \colon \mathbb{M} \m{od}(R) & \to & \mathbb{M} \m{od}(R') \\
  \M & \rightsquigarrow & \M \otimes_R R'.
\end{eqnarray*}
For an $R$-category $\M$, we put $\M_{R'} \coloneqq \M \otimes_R R'$ for short. Similarly for an object $M$ of $\M$, we denote by $M_{R'}$ the object $M \otimes_R R'$ of $\M_{R'}$ corresponding to $M$. Then $\M_{R'}$ is the category of objects in $\M$ and the homsets
\begin{eqnarray*}
  \m{Hom}_{\M_{R'}}(M^1_{R'},M^2_{R'}) \cong_{R'} \m{Hom}_{\M}(M^1,M^2) \otimes_R R'
\end{eqnarray*}
for objects $M^1$ and $M^2$ of $\M$. The natural $R$-linear functor $\m{id}_{\M} \otimes_R R' \colon \M \to \M_{R'}$ has the universal property that for any $R'$-category $\M'$ and any $R$-linear functor $\Phi \colon \M \to \M'$, there exists an $R'$-linear functor $\Phi_{R'} \colon \M_{R'} \to \M'$ unique up to $R'$-linear natural equivalences such that $\Phi_{R'} \circ (\m{id}_{\M} \otimes_R R')$ is $R$-linearly naturally equivalent to $\Phi$. In particular, for a $k^{\circ}$-category $\M$, the $k$-category $\M_k$ is the quotient category of $\M$ by the Serre subcategory of objects admitting no $k^{\circ}$-torsionfree endomorphisms.

\subsection{Banach Spaces and Bounded Homomorphisms}
\label{Banach Spaces and Bounded Homomorphisms}

\begin{prp}
\label{int ban - ban}
Suppose that the valuation of $k$ is nontrivial and discrete. The covariant functor
\begin{eqnarray*}
  \m{BanC}^{\m{crt}}(k)_k & \to & \m{Ban}(k) \\
  V_k & \rightsquigarrow & V
\end{eqnarray*}
is an equivalence of categories.
\end{prp}

\begin{proof}
The assertion is verified in the proof of \cite{ST} Theorem 1.2 when $k$ is a local field of characteristic $0$, and the proof is valid for an arbitrary complete nontrivial discrete valuation field $k$.
\end{proof}

\begin{dfn}
A unitarisable Banach $k$-linear representation of $G$ is a pair $(V,\pi)$ of a Banach $k$-vector space $V$ and a continuous map $\pi \colon G \times V \to V$ such that $(V,\pi)$ is a unitary Banach $k$-linear representation of $G$ with respect to some equivalent norm on $V$. Denote by $\m{Rep}_{G}^{\m{unit}}(\m{Ban}(k))$ the category of unitarisable Banach $k$-linear representations of $G$ and continuous $G$-equivariant $k$-linear homomorphisms.
\end{dfn}

\begin{prp}
\label{rep int ban - ban}
Suppose that the valuation of $k$ is nontrivial and discrete. The covariant functor
\begin{eqnarray*}
  \m{Rep}_{G} \left( \m{BanC}^{\m{crt}}(k) \right)_k & \to & \m{Rep}_{G}^{\m{unit}}(\m{Ban}(k)) \\
  (V,\pi)_k & \rightsquigarrow & (V_k,\pi_k)
\end{eqnarray*}
is an equivalence of categories.
\end{prp}

\begin{proof}
Denote the functor by $\iota$. The forgetful functors $\m{Rep}_{G}(\m{BanC}^{\m{crt}}(k)) \to \m{BanC}^{\m{crt}}(k)$ and $\m{Rep}_{G}^{\m{unit}}(\m{Ban}(k)) \to \m{Ban}(k)$ are faithful $k^{\circ}$-linear functors, and so is the induced functor $\m{Rep}_{G}(\m{BanC}^{\m{crt}}(k))_k \to \m{BanC}^{\m{crt}}(k)_k$ because $k$ is a flat $k^{\circ}$-module. The functor $\iota$ is compatible with the equivalence $\m{BanC}^{\m{crt}}(k)_k \to \m{Ban}(k)$ in Proposition \ref{int ban - ban}, and hence is faithful. For the fullness, let $(V^1,\pi^1)$ and $(V^2,\pi^2)$ be objects of $\m{Rep}_{G}(\m{BanC}(K))$. For a morphism $A \colon (V^1_k,\pi^1_k) \to (V^2_k,\pi^2_k)$ in $\m{Rep}_{G}^{\m{unit}}(\m{Ban}(k))$, since the underlying homomorphism $A \colon V^1_k \to V^2_k$ is continuous, there are a contractive homomorphism $B \colon V^1 \to V^2$ and an $a \in k^{\times}$ such that $a B_k = A$ in $\m{Ban}(k)$. The actions of $G$ are $k$-linear, and hence $B$ is $G$-equivariant. Therefore $B$ gives a morphism in $\m{Rep}_{G}(\m{BanC}^{\m{crt}}(k))$ with $\iota(a B_k) = A$ in $\m{Rep}_{G}^{\m{unit}}(\m{Ban}(k))$. Therefore $\iota$ is full.

Let $(\V,\varpi)$ be an object of $\m{Rep}_{G}^{\m{unit}}(\m{Ban}(k))$. We denote by $\n{\cdot}$ the norm of $\V$. Since $(\V,\varpi)$ is unitarisable, there is an equivalent norm $\n{\cdot}'$ on $\V$ preserved by the action of $G$. Take $R_1 > R_2 > 0$ with $R_2 \n{\cot}' \leq \n{\cdot} \leq R_1 \n{cdot}'$. Then for any $(g,v) \in G \times \V$, the equality $\n{gv}' = \n{v}'$ implies $R_1^{-1}R_2 \n{v} \leq \n{gv} \leq R_1R_2^{-1} \n{v}$. Therefore the orbit $Gv \subset \V$ is bounded. Set $\n{v}'' \coloneqq \sup_{g \in G} \n{gv}$. For any $v_1, v_2 \in \V$, the equality
\begin{eqnarray*}
  \n{v_1 + v_2}'' & = & \sup_{g \in G} \n{gv_1 + gv_2} \leq \sup_{g \in G} \max \Ens{\n{gv_1},\n{gv_2}} \\
  & = & \max \Ens{\sup_{g \in G} \n{gv_1}, \sup_{g \in G} \n{gv_2}} = \max \Ens{\n{v_1}'', \n{v_2}''}
\end{eqnarray*}
holds, and for any $(a,v) \in k \times \V$, the equality
\begin{eqnarray*}
  \n{av}'' = \sup_{g \in G} \n{g(av)} = \sup_{g \in G} \v{a} \ \n{gv} = \v{a} \sup_{g \in G} \n{gv} = \v{a} \ \n{v}''
\end{eqnarray*}
holds. Therefore $\n{c\dot}''$ is a seminorm on the underlying $k$-vector space of $\V$. Moreover, for any $v \in \V$, the inequality $R_1^{-1}R_2 \n{v} \leq \n{gv} \leq R_1R_2^{-1} \n{v}$ for every $g \in G$ guarantees $R_1^{-1}R_2 \n{v} \leq \n{v}'' \leq R_1R_2^{-1} \n{v}$. Thus $\n{\cdot}''$ is a norm equivalent to $\n{\cdot}$. In particular, $\n{\cdot}''$ is complete. The map
\begin{eqnarray*}
  \n{\cdot}_{\V,\varpi} \colon \V & \to & [0, \infty ) \\
  v & \mapsto & \inf \Set{\v{a}}{a \in k, \v{a} \geq \n{v}}.
\end{eqnarray*}
is also a complete norm of the underlying $k$-vector space of $\V$ equivalent to $\n{\cdot}''$ by a similar calculation. Since the valuation of $k$ is discrete, the image of $\n{\cdot}_{\V,\varpi}$ is contained in $\v{k}$. Denote by $V$ the Banach $k$-vector space whose underlying $k$-vector space is that of $\V$ and whose norm is $\n{\cdot}_{\V,\varpi}$. We set $\pi \coloneqq \varpi$. The strictly Cartesian unitary Banach $k$-linear representation $(V,\pi)$ satisfies $(V_k,\pi_k) \cong (\V,\varpi)$ in $\m{Rep}_{G}^{\m{unit}}(\m{Ban}(k))$. The construction of $(V,\pi)$ is natural on $(\V,\varpi)$, and since $\iota$ is fully faithful, it gives the inverse functor
\begin{eqnarray*}
  \m{Rep}_{G}^{\m{unit}}(\m{Ban}(k)) \to \m{Rep}_{G} \left( \m{BanC}^{\m{crt}}(k) \right)_k
\end{eqnarray*}
of $\iota$. We conclude that $\iota$ is an equivalence.
\end{proof}

\begin{dfn}
An integralisable Banach $k[[G]]$-module is a $k[[G]]$-module $V$ endowed with a complete norm $\n{\cdot}$ as a $k$-vector space such that $V$ is an integral Banach $k[[G]]$-module with respect to some norm $\n{\cdot}'$ equivalent to $\n{\cdot}$. Denote by $\m{Ban}^{\m{int}}(k[[G]])$ the category of integralisable Banach $k[[G]]$-modules and continuous $k[[G]]$-linear homomorphisms.
\end{dfn}

\begin{lmm}
\label{equicontinuity}
Let $V$ be a $k[[G]]$-module endowed with a complete norm $\n{\cdot}$. Then $V$ is an integralisable Banach $k[[G]]$-module if and only if $V$ is a topological $k^{\circ}[[G]]$-module with respect to the norm topology and there is an $R > 0$ such that $\n{\mu v} \leq R \n{v}$ for any $(\mu,v) \in k^{\circ}[[G]] \times V$.
\end{lmm}

\begin{proof}
Suppose that $V$ is an integralisable Banach $k[[G]]$-module. Take a norm $\n{\cdot}'$ equivalent to $\n{\cdot}$ for which $V$ is an integral Banach $k[[G]]$-module. Then $V$ is a topological $k^{\circ}[[G]]$-module with respect to the norm topology. Since $\n{\cdot}'$ is equivalent to $\n{\cdot}$, there is an $R > 0$ such that $\n{\cdot}' \leq R \n{\cdot}$. Every element of $k^{\circ}[[G]]$ corresponds to a contractive endomorphism of $V$ with respect to $\n{\cdot}'$ by the scalar multiplication, and hence $\n{\mu v} \leq R \n{\mu v} \leq R \n{v}$ for any $(\mu,v) \in k^{\circ}[[G]] \times V$.

Suppose that $V$ is a topological $k^{\circ}[[G]]$-module and there is an $R > 0$ such that $\n{\mu v} \leq R \n{v}$ for any $(\mu,v) \in k^{\circ}[[G]] \times V$. For each $v \in V$, set $\n{v}' \coloneqq \sup_{\mu \in k^{\circ}[[G]]} \n{\mu v}$. Then $\n{\cdot}' \colon V \to [0, \infty )$ is a norm of a $k$-vector space. It is easy to see that $\n{\cdot} \leq \n{\cdot}' \leq R \n{\cdot}$, and hence $\n{\cdot}'$ is equivalent to $\n{\cdot}$. Every element of $k^{\circ}[[G]]$ corresponds to a contractive endomorphism of $V$ with respect to $\n{\cdot}'$ by the scalar multiplication, and thus $V$ is an integral Banach $k[[G]]$-module with respect to $\n{\cdot}'$.
\end{proof}

\begin{prp}
\label{int iwa ban - iwa ban}
Suppose that the valuation of $k$ is nontrivial and discrete. The covariant functor
\begin{eqnarray*}
  \m{BanC}^{\m{crt}}(k^{\circ}[[G]])_k & \to & \m{Ban}^{\m{int}}(k[[G]]) \\
  (V,\pi)_k & \rightsquigarrow & (V_k,\pi_k)
\end{eqnarray*}
is an equivalence of categories.
\end{prp}

\begin{proof}
The assertion follows from Proposition \ref{int ban - ban} and Lemma \ref{equicontinuity} by a similar argument as in Proposition \ref{rep int ban - ban}.
\end{proof}

\begin{prp}
Suppose that the valuation of $k$ is nontrivial and discrete. The covariant functor
\begin{eqnarray*}
  (\cdot)|_G \colon \m{Ban}^{\m{int}}(k[[G]]) & \to & \m{Rep}_{G}^{\m{unit}}(\m{Ban}(k)) \\
  (V,\pi) & \rightsquigarrow & (V,\pi|_G)
\end{eqnarray*}
is an equivalence of categories.
\end{prp}

\begin{proof}
The functor is the composite of the equivalences in Corollary \ref{int iwa - rep ban 4}, Proposition \ref{rep int ban - ban}, and Proposition \ref{int iwa ban - iwa ban}, and hence is an equivalence.
\end{proof}

\subsection{Compactly Generated Locally Convex Spaces}
\label{Compactly Generated Locally Convex Spaces}

\begin{dfn}
\label{compactly generated}
A compactly generated locally convex $k$-vector space is a pair $(W_{\tau_W},\n{\cdot})$ of a Hausdorff locally convex $k$-vector space $W_{\tau_W} = (W,\tau_W)$ and a norm $\n{\cdot}$ on $W$ such that the closed unit ball $M \subset W$ with respect to $\n{\cdot}$ is compact with respect to $\tau_W|_M$, and such that $\tau_W$ is the strongest locally convex topology on $W$ for which the embedding $M_{\tau_W|_M} \hookrightarrow W$ is continuous. We denote by $\m{TMod}_{\m{lc,cg}}(k)$ the category of compactly generated locally convex $k$-vector spaces and $k$-linear homomorphisms continuous with respect to the locally convex topologies.
\end{dfn}

\begin{prp}
For a compactly generated locally convex $k$-vector space $(W_{\tau_W},\n{\cdot})$, the norm topology of $W$ is stronger than $\tau_W$.
\end{prp}

\begin{proof}
Let $U \subset W_{\tau_W}$ be an open neighbourhood of $0$. Denote by $M \subset W$ the closed unit ball with respect to $\n{\cdot}$. For any $w \in W$, the scalar multiplication $k \to W \colon a \mapsto aw$ is continuous. It implies $\bigcup_{a \in k^{\times}} aU = W$. Since $M$ is compact, there is an $a \in k^{\times}$ such that $M \subset aU$. Thus $U$ contains the closed ball centred at $0$ and of radius $\v{a}^{-1}$. We conclude that the norm topology of $W$ is stronger than $\tau_W$.
\end{proof}

For a compact Hausdorff flat linear topological $k^{\circ}$-module $M_{\tau_M}$, we endow $M \otimes_{k^{\circ}} k$ with the strongest locally convex topology $\tau^M$ for which the embedding $M_{\tau_M} \hookrightarrow M \otimes_{k^{\circ}} k$ is continuous. By the argument before Lemma 1.4 in \cite{ST}, $\tau^M$ is Hausdorff, and hence $\tau_M$ coincides with $\tau^M|_M$. For each $m \in M \otimes_{k^{\circ}} k$, we set $\n{m}_M \coloneqq \inf \set{\v{a}}{a \in k^{\times}, a^{-1}m \in M}$. Since $M$ is adically separated by Lemma \ref{compact - adic}, $\n{m}_M = 0$ if and only if $m = 0$. If the valuation of $k$ is discrete, then the closed unit ball of $M \otimes_{k^{\circ}} k$ with respect to $\n{\cdot}_M$ coincides with $M$. Thus $((M \otimes_{k^{\circ}} k)_{\tau^M},\n{\cdot}_M)$ is a compactly generated locally convex $k$-vector space. The natural correspondence $M \rightsquigarrow ((M \otimes_{k^{\circ}} k)_{\tau^M},\n{\cdot}_M)$ gives a $k^{\circ}$-linear functor $\m{TMod}_{\m{flat,lin}}^{\m{cpt,sep}}(k^{\circ}) \to \m{TMod}_{\m{lc,cg}}(k)$.

\begin{prp}
\label{cpt lin - cg lc}
Suppose that $k$ is a local field. The covariant functor
\begin{eqnarray*}
  \m{TMod}_{\m{flat,lin}}^{\m{cpt,sep}}(k^{\circ})_k & \to & \m{TMod}_{\m{lc,cg}}(k) \\
  M_k & \rightsquigarrow & ((M \otimes_{k^{\circ}} k)_{\tau^M},\n{\cdot}_M)
\end{eqnarray*}
is an equivalence of categories.
\end{prp}

\begin{proof}
Denote the functor by $\iota$. The $k^{\circ}$-linear functor $\m{TMod}_{\m{flat,lin}}^{\m{cpt,sep}}(k^{\circ}) \to \m{TMod}_{\m{lc,cg}}(k)$ is faithful because $M \otimes_{k^{\circ}} k$ is generated by $M$ as a $k$-vector space for any object $M$ of $\m{TMod}_{\m{flat,lin}}^{\m{cpt,sep}}(k^{\circ})$, and hence so is $\iota$ because $k$ is a flat $k^{\circ}$-module. The fullness follows from Lemma 1.5 (iii) in \cite{ST}. We remark that the proof of Lemma 1.5 is valid for a local field of any characteristic.

Let $(W_{\tau_W},\n{\cdot})$ be an object of $\m{TMod}_{\m{lc,cg}}(k)$. Denote by $M \subset W$ the closed unit ball with respect to $\n{\cdot}$ endowed with $\tau_W|_M$. By the definition of $\tau^M$, $(M \otimes_{k^{\circ}} k)_{\tau^M}$ coincides with $W_{\tau_W}$. Since the valuation of $k$ is discrete, it is easy to see that $\n{\cdot}_M$ is equivalent to $\n{\cdot}$. Thus $((M \otimes_{k^{\circ}} k)_{\tau^M},\n{\cdot}_M)$ is naturally isomorphic to $(W_{\tau_W},\n{\cdot})$. The construction of $M$ is natural, and since $\iota$ is fully faithful, it gives the inverse functor
\begin{eqnarray*}
  \m{TMod}_{\m{lc,cg}}(k) \to \m{TMod}_{\m{flat,lin}}^{\m{cpt,sep}}(k^{\circ})_k
\end{eqnarray*}
of $\iota$. We conclude that $\iota$ is an equivalence.
\end{proof}

\begin{prp}
\label{replacement of the norm}
Suppose that the valuation of $k$ is nontrivial and discrete. For an object $(W_{\tau_W},\n{\cdot})$ of $\m{TMod}_{\m{lc,cg}}(k)$ and a norm $\n{\cdot}'$ equivalent to $\n{\cdot}$, let $M, M' \subset W$ be the closed unit balls with respect to $\n{\cdot}$ and $\n{\cdot}'$ respectively. Then $M'|_{\tau_W|_{M'}}$ is compact, and $\tau_W$ coincides with the strongest locally convex topology on $W$ for which the embedding $M'_{\tau_W|_{M'}} \hookrightarrow W$ is continuous.
\end{prp}

\begin{proof}
Denote by $\tau'_W$ the strongest locally convex topology on $W$ for which the embedding $M'_{\tau_W|_{M'}} \hookrightarrow W$ is continuous. Take $a_1,a_2 \in k^{\circ}$ with $a_1 \n{\cdot} \leq \n{\cdot}' \leq a_2 \n{\cdot}$. Then $a_2^{-1} M \subset M' \subset a_1^{-1} M$. Since $M$ is compact and Hausdorff, so is $a_1^{-1}M$. By the equivalence of $\n{\cdot}$ and $\n{\cdot}'$, $\n{\cdot}'$ defines the topology of $a_1^{-1}M$. Therefore $M' \subset a_1^{-1}M$ is clopen, and is compact. Since the embedding $M'_{\tau'_W|_{M'}} \hookrightarrow W_{\tau_W}$ is continuous, $\tau'_W$ is stronger than or equal to $\tau_W$ by the universality of $\tau'_W$. Moreover, the embedding $M_{\tau_W} \hookrightarrow W_{\tau'_W}$ coincides with the composite
\begin{eqnarray*}
  M_{\tau_W} \stackrel{a_2^{-1}}{\to} a_2^{-1} M_{\tau_W} \hookrightarrow M'_{\tau'_W} \hookrightarrow W_{\tau'_W} \stackrel{a_2}{\to} W_{\tau'_W},
\end{eqnarray*}
and hence is continuous. Thus $\tau_W = \tau'_W$ by the universality of $\tau_W$.
\end{proof}

Note that the functor forgetting the norm from compactly generated locally convex spaces is fully faithful and essentially surjective, but is not an equivalence of category because we do not fix a universe so that the categories are small. This is why we equip such a locally convex space with a non-canonical norm.

\begin{dfn}
A unitarisable compactly generated locally convex $k$-linear representation of $G$ is a pair $((W_{\tau_W},\n{\cdot}),\rho)$ of a compactly generated locally convex $k$-vector space $(W_{\tau_W},\n{\cdot})$ and a continuous $k$-linear action $\rho \colon G \times W_{\tau_W} \to W_{\tau_W}$ preserving some norm equivalent to $\n{\cdot}$. Denote by $\m{Rep}_{G}^{\m{unit}}(\m{TMod}_{\m{lc,cg}}(k))$ the category of unitarisable compactly generated locally convex $k$-linear representations of $G$ and $G$-equivariant homomorphisms of compactly generated locally convex $k$-vector spaces.
\end{dfn}

\begin{lmm}
\label{rep int - rep lc}
Suppose that $k$ is a local field. The covariant functor
\begin{eqnarray*}
  \m{Rep}_{G}(\m{TMod}_{\m{flat,lin}}^{\m{cpt,sep}}(k^{\circ}))_k & \to & \m{Rep}_{G}^{\m{unit}}(\m{TMod}_{\m{lc,cg}}(k)) \\
  (M,\rho)_k & \rightsquigarrow & (((M \otimes_{k^{\circ}} k)_{\tau^M},\n{\cdot}_M),\rho \otimes_{k^{\circ}} k)
\end{eqnarray*}
is an equivalence of categories.
\end{lmm}

\begin{proof}
Denote by $\iota$ the functor. The full faithfulness of $\iota$ is verified observing the forgetful functors in a usual way by Proposition \ref{cpt lin - cg lc} and \cite{ST} Lemma 1.5 (iii). Note that $k$ is assumed to be of characteristic $0$ in \cite{ST} Lemma 1.5 but the proof is valid for positive characteristic. Let $((W_{\tau_W},\n{\cdot}),\rho_W)$ be an object of $\m{Rep}_{G}^{\m{unit}}(\m{TMod}_{\m{lc,cg}}(k))$. In a similar way as in the proof of Proposition \ref{rep int ban - ban}, there is a natural construction of a norm $\n{\cdot}_{\rho_W}$ equivalent to $\n{\cdot}$ preserved by the action of $G$. Denote by $M' \subset W$ the closed unit ball with respect to $\n{\cdot}_{\rho_W}$ endowed with $\tau_{\rho_W}|_{M'}$. Then $(((M' \otimes_{k^{\circ}} k)_{\tau^{M'}},\n{\cdot}_{M'}),\rho \otimes_{k^{\circ}} k)$ coincides with $((W_{\tau_W},\n{\cdot}_{\rho_W}),\rho_W)$ by Proposition \ref{replacement of the norm}, and is isomorphic to $((W_{\tau_W},\n{\cdot}),\rho_W)$ in $\m{Rep}_{G}^{\m{unit}}(\m{TMod}_{\m{lc,cg}}(k))$. The construction of $(M,\rho)$ is natural on $((W_{\tau_W},\n{\cdot}),\rho_W)$, and hence by the full faithfulness of $\iota$, it gives a contravariant functor
\begin{eqnarray*}
  \m{Rep}_{G}^{\m{unit}}(\m{TMod}_{\m{lc,cg}}(k)) \to \m{Rep}_{G}(\m{TMod}_{\m{flat,lin}}^{\m{cpt,sep}}(k^{\circ}))_k
\end{eqnarray*}
which is the inverse of $\iota$. We conclude that $\iota$ is an equivalence.
\end{proof}

\begin{dfn}
\label{integralisable}
An integralisable compactly generated locally convex $k[[G]]$-module is a $k[[G]]$-module $W$ endowed with a locally convex topology $\tau_W$ and a complete norm $\n{\cdot}$ as a $k$-vector space such that $W$ is a compactly generated locally convex $k$-vector space with with respect to $\tau_W$ and $\n{\cdot}$, and such that the closed unit ball of $W$ with respect to some norm equivalent to $\n{\cdot}$ is a $k^{\circ}[[G]]$-submodule. Denote by $\m{TMod}_{\m{lc,cg}}^{\m{int}}(k[[G]])$ the category of integralisable compactly generated locally convex $k[[G]]$-modules and continuous $k[[G]]$-linear homomorphisms.
\end{dfn}

\begin{prp}
\label{int iwa - lc iwa}
Suppose that $k$ is a local field. The covariant functor
\begin{eqnarray*}
  \m{TMod}_{\m{flat,lin}}^{\m{cpt,sep}}(k^{\circ}[[G]])_k & \to & \m{TMod}_{\m{lc,cg}}^{\m{int}}(k[[G]]) \\
  M_k & \rightsquigarrow & ((M \otimes_{k^{\circ}} k)_{\tau^M},\n{\cdot}_M)
\end{eqnarray*}
is an equivalence of categories.
\end{prp}

\begin{proof}
The inverse is constructed in a similar way as in Lemma \ref{rep int - rep lc} using  Proposition \ref{cpt lin - cg lc}, Proposition \ref{replacement of the norm}, and \cite{ST} Lemma 1.5 (iii).
\end{proof}

Thus we identify $M_k$ and $((M \otimes_{k^{\circ}} k)_{\tau^M},\n{\cdot}_M)$ through the equivalence. We remark that the norm of an object of $\m{TMod}_{\m{lc,cg}}^{\m{int}}(k[[G]])$ is not an invariant of an isomorphism class, and hence $M_k$ does not possess the information of the canonical lattice $M \subset M_k$. It just restore the isomorphism class of the lattice.

\begin{prp}
Suppose that $k$ is a local field. The covariant functor
\begin{eqnarray*}
  (\cdot)|_G \colon \m{TMod}_{\m{lc,cg}}^{\m{int}}(k[[G]]) & \to & \m{Rep}_{G}^{\m{unit}}(\m{TMod}_{\m{lc,cg}}(k)) \\
  (W,\rho) & \rightsquigarrow & (W,\rho|_G)
\end{eqnarray*}
is an equivalence of categories.
\end{prp}

\begin{proof}
The functor is the composite of the equivalences in Proposition \ref{int iwa - rep cpt 3}, Lemma \ref{rep int - rep lc}, and Proposition \ref{int iwa - lc iwa}.
\end{proof}

\subsection{Linear-Type Duality}
\label{Linear-Type Duality}

\begin{thm}
\label{iwa - rep ban}
Suppose that $k$ is a local field and $G$ is a locally profinite group. There is a contravariant equivalence
\begin{eqnarray*}
  \mathbb{D} \colon \m{TMod}_{\m{lc,cg}}^{\m{int}}(k[[G]]) & \to & \m{Rep}_{G}^{\m{unit}}(\m{Ban}(k))
\end{eqnarray*}
of $\m{Mod}(k)$-enriched categories.
\end{thm}

\begin{proof}
The equivalence is given as the composite of the equivalences in Theorem \ref{int iwa - rep ban}, Proposition \ref{rep int ban - ban}, and Proposition \ref{int iwa - lc iwa}.
\end{proof}

For a unitarisable Banach $k$-linear representation $(V,\pi)$ of $G$, we denote by $V^{\mathbb{D}}$ the underlying compactly generated locally convex $k$-vector space of $(V,\pi)^{\mathbb{D}}$, and by $\pi^{\mathbb{D}}$ the scalar multiplication $k[[G]] \times V^{\mathbb{D}} \to V^{\mathbb{D}}$. A Banach $k$-linear representation of $G$ is said to be irreducible if it contains no nontrivial closed $G$-stable $k$-vector subspace. Similarly, an integral compactly generated locally convex $k[[G]]$-module is said to be simple if it contains no nontrivial closed $k[[G]]$-submodule. The duality theorem yields a criterion for the irreducibility of a unitaralisable Banach $k$-linear representation.

\begin{crl}
\label{irr - smp}
Suppose that $k$ is a local field and $G$ is a locally profinite group. Then a unitarisable Banach $k$-linear representation $(V,\pi)$ is irreducible if and only if $(V,\pi)^{\mathbb{D}}$ is a simple integral compactly generated locally convex $k[[G]]$-module.
\end{crl}

\begin{proof}
Suppose that $(V,\pi)^{\mathbb{D}}$ is simple. Let $V_1 \subset V$ be a proper closed $G$-stable $k$-vector subspace. By Hahn--Banach theorem (\cite{ST} Proposition 9.2), the homomorphism $V^{\mathbb{D}} \to V_1^{\mathbb{D}}$ induced by the inclusion $V_1 \hookrightarrow V$is surjective and the kernel is a non-zero closed $k[[G]]$-submodule of $V^{\mathbb{D}}$. Therefore $V_1^{\mathbb{D}} = 0$ and $V_1 = 0$ by the equivalence of $\mathbb{D}$. Thus $(V,\pi)$ is irreducible.

Suppose that $(V,\pi)$ is irreducible. Denote by $\tau_{V^{\mathbb{D}}}$ the locally convex topology of $V^{\mathbb{D}}$. Let $W \subset V^{\mathbb{D}}$ be a proper closed $k[[G]]$-submodule. Since $(V,\pi)^{\mathbb{D}}$ is integralisable, there is a norm $\n{\cdot}'$ equivalent to the underlying norm of $V^{\mathbb{D}}$ such that the closed unit ball $M' \subset V^{\mathbb{D}}$ with respect to $\n{\cdot}'$ is a compact $k^{\circ}[[G]]$-submodule. By Proposition \ref{replacement of the norm}, $\tau_{V^{\mathbb{D}}}$ coincides with $\tau^{M'}$. Put $M'' \coloneqq W \cap M'$. Since $W \subset V^{\mathbb{D}}$ is closed, $M'/M''$ is a compact Hausdorff $k^{\circ}$-flat $k^{\circ}$-linear topological $k^{\circ}[[G]]$-module with respect to the quotient topology $\m{quo}$. Note that a $k^{\circ}$-module is flat if and only if it is torsionfree. Therefore $((V^{\mathbb{D}}/W)_{\tau^{M'/M''}},\n{\cdot}_{M'/M''})$ is an integralisable compactly generated locally convex $k$-vector space. Since the canonical projection $M'|_{\tau_{V^{\mathbb{D}}}|_{M'}} \twoheadrightarrow (M'/M'')_{\m{quo}}$ is continuous, so is $V^{\mathbb{D}}_{\tau_{V^{\mathbb{D}}}} \twoheadrightarrow (V^{\mathbb{D}}/W)_{\tau^{M'/M''}}$ by the universality of $\tau^{M'} = \tau_{V^{\mathbb{D}}}$. By Theorem \ref{iwa - rep ban}, the complex $W_{\tau_{V^{\mathbb{D}}}|_W} \hookrightarrow V^{\mathbb{D}}_{\tau_{V^{\mathbb{D}}}} \twoheadrightarrow (V^{\mathbb{D}}/W)_{\tau^{M'/M''}}$ is obtained as the dual of a complex $(V_1,\pi_1) \to (V,\pi) \to (V_2,\pi_2)$ for some unitarisable Banach $k$-linear representations $(V_1,\pi_1)$ and $(V_2,\pi_2)$. The canonical projection $V^{\mathbb{D}} \twoheadrightarrow V^{\mathbb{D}}/W$ is a nonzero homomorphism, and hence so is the dual $V_1 \to V$. Since $(V,\pi)$ is irreducible, the image of $V_1 \to V$ is dense. Then its dual $V^{\mathbb{D}} \to V^{\mathbb{D}}/W$ is injective, and hence $W = O$. Thus $(V,\pi)^{\mathbb{D}}$ is simple.
\end{proof}

\section{Applications}
\label{Applications}

Suppose that $k$ is a local field and $G$ is a locally profinite group. Let $H \subset G$ be a closed subgroup. In this section, we give an interpretation of operations on unitary Banach representations as those on compact Hausdorff $k^{\circ}$-flat $k^{\circ}$-linear topological Iwasawa modules. As a consequence, we give an explicit description of a continuous parabolic induction of a unitary Banach representation.

\subsection{Several Operations}
\label{Several Operations}

Let $V$ be a unitary (resp.\ unitarisable) Banach $k$-linear representation, and $M$ (resp.\ $W$) a compact Hausdorff $k^{\circ}$-flat $k^{\circ}$-linear topological $k^{\circ}[[G]]$-module (resp.\ an integralisable compactly generated locally convex $k[[G]]$-module) corresponding to each other. As we have already verified in the proof of Corollary \ref{irr - smp},  there is an anti-order-preserving one-to-one corresponding between closed subrepresentations $V_1 \subset V$ and closed $k^{\circ}[[G]]$-submodules $M_1 \subset M$ (resp.\ closed $k[[G]]$-submodules $W_1 \subset W$) given by the duality functors in Theorem \ref{int iwa - rep ban} and Theorem \ref{iwa - rep ban}. This interprets closed subrepresentations and quotient representations as quotient submodules and closed submodules respectively.

\vspace{0.1in}
For a Banach $k$-linear representation $V$ of $G$, we denote by $\m{Res}_{H}^{G}(V)$ the Banach representation obtained by restricting the action $G \times V \to V$ to $H \times V \to V$. This gives contravariant functors
\begin{eqnarray*}
  & & \m{Res}_{H}^{G} \colon \m{Rep}_{G}(\m{BanC}(k)) \to \m{Rep}_{H}(\m{BanC}(k)) \\
  & & \m{Res}_{H}^{G} \colon \m{Rep}_{G}(\m{BanC}^{\m{crt}}(k)) \to \m{Rep}_{H}(\m{BanC}^{\m{crt}}(k)) \\
  & & \m{Res}_{H}^{G} \colon \m{Rep}_{G}^{\m{unit}}(\m{Ban}(k)) \to \m{Rep}_{H}^{\m{unit}}(\m{Ban}(k)),
\end{eqnarray*}
and their duals are the covariant functors
\begin{eqnarray*}
  & & \m{Res}_{k^{\circ}[[H]]}^{k^{\circ}[[G]]} \colon \m{TMod}_{\m{flat,lin}}^{\m{cpt,sep}}(k^{\circ}[[G]]) \to \m{TMod}_{\m{flat,lin}}^{\m{cpt,sep}}(k^{\circ}[[H]]) \\
  & & \m{Res}_{k[[H]]}^{k[[G]]} \colon \m{TMod}_{\m{lc,cg}}^{\m{int}}(k[[G]]) \to \m{TMod}_{\m{lc,cg}}^{\m{int}}(k[[H]])
\end{eqnarray*}
given by restriction of scalar.

\vspace{0.1in}
Suppose that $H \backslash G$ is compact with respect to the quotient topology. Then there is a compact open representative $K \subset G$ of the canonical projection $G \twoheadrightarrow H \backslash G$. Indeed, take an compact open subgroup $G_0 \subset G$. Since $G \twoheadrightarrow H \backslash G$ is an open map, there are an $n \in \N$ and $g_1, \ldots, g_n \in G$ such that $H \backslash G = \bigcup_{i = 1}^{n} HG_0 g_i$. It implies that $\bigcup_{i = 1}^{n} G_0 g_i$ is a compact open representative. For a unitarisable Banach $k$-linear representation $V_0$ of $H$, we set
\begin{eqnarray*}
  \m{Ind}_{H}^{G}(V_0) \coloneqq \Set{f \in \m{C}_{\m{bd}}(G,V_0)}{f(hv) = hf(v), {}^{\forall}(h,v) \in H \times V_0}.
\end{eqnarray*}
Note that for any continuous map $f \colon G \to V_0$ with $f(hv) = hf(v)$ for any $(h,v) \in H \times V_0$, $f$ is bounded. Indeed, take a norm $\n{\cdot}'$ equivalent to the original norm $\n{\cdot}$ of $V_0$ such that $V_0$ is a unitary representation of $H$ with respect to $\n{\cdot}'$. For any $(h,g) \in H \times G$, we have $\n{f(hg)} = \n{hf(g)}' = \n{f(g)}'$, and hence the continuous map $\n{f(\cdot)}' \colon G \to [0, \infty ) \colon g \mapsto \n{f(g)}'$ is left $H$-invariant. Therefore it induces a continuous map $H \backslash G \to [0, \infty )$ by the universality of the quotient topology. Since $H \backslash G$ is compact, the image is bounded. By the equivalence between $\n{\cdot}$ and $\n{\cdot}'$, $f$ is bounded. Now we consider the isometric action $G \times \m{C}_{\m{bd}}(G,V_0) \to \m{C}_{\m{bd}}(G,V_0)$ given by setting $(gf)(g') \coloneqq f(g'g)$ for each $(f,g,g') \in \m{C}_{\m{bd}}(G,V_0) \times G \times G$. Beware that the action is not strongly continuous when $G$ is not compact. Then $\m{Ind}_{H}^{G}(V_0) \subset \m{C}_{\m{bd}}(G,V_0)$ is a closed $G$-equivariant $k$-vector subspace. Moreover, the restriction $G \times \m{Ind}_{H}^{G}(V_0) \to \m{Ind}_{H}^{G}(V_0)$ is strongly continuous and hence $\m{Ind}_{H}^{G}(V_0)$ is a unitary Banach $k$-linear representation of $G$. Indeed, take a compact open representative $K \subset G$ of the canonical projection $G \twoheadrightarrow H \backslash G$. Then the norm of $\m{Ind}_{H}^{G}(V_0)$ is equivalent to the supremum norm $\n{\cdot}_K$ on $K$ by the unitarisability of the action of $H$ on $V_0$. Since $K \subset G$ is a compact open subspace, there is a compact open subgroup $G_0 \subset G$ with $K G_0 = K$. Then the restricted action $G_0 \times \m{Ind}_{H}^{G}(V_0) \to \m{Ind}_{H}^{G}(V_0)$ preserves $\n{\cdot}_K$. Since a continuous function on a compact uniform space is uniformly continuous, the action of $G_0$ is separately continuous. The family of isometric operators with respect to $\n{\cdot}_K$ is equicontinuous by Banach--Steinhaus theorem (\cite{Sch} Corollary 6.16), and hence the action $G \times \m{Ind}_{H}^{G}(V_0) \to \m{Ind}_{H}^{G}(V_0)$ is strongly continuous.

\vspace{0.1in}
For an object $M_0$ of $\m{TMod}_{\m{flat,lin}}^{\m{cpt,sep}}(k^{\circ}[[H]])$, we express $\m{Ind}_{H}^{G}(M_0^D)^D$ explicitly by $G$ and $M_0$. By a similar argument as in Lemma \ref{rep - unf}, giving a continuous map $G \to M_0^D$ is equivalent to giving a continuous map $G \times M_0 \to k$ which is $k^{\circ}$-linear with respect to the variable in $M_0$. Therefore we identify $\m{C}_{\m{bd}}(G,M_0^D)$ as a closed $k$-vector subspace of $\m{C}_{\m{bd}}(G \times M_0,k)$. The inclusion $\m{Ind}_{H}^{G}(M_0^D) \hookrightarrow \m{C}_{\m{bd}}(G \times M_0,k)$ is $G$-equivariant with respect to the isometric action $G \times \m{C}_{\m{bd}}(G \times M_0,k) \to \m{C}_{\m{bd}}(G \times M_0,k)$ given by setting $(gf)(g',m) \coloneqq f(g'g,m)$ for each $(g,f,g',m) \in G \times \m{C}_{\m{bd}}(G \times M_0,k) \times G \times M_0$. The action is not strongly continuous when $G$ is not compact. Through the inclusion, $\m{Ind}_{H}^{G}(M_0^D)$ corresponds to the closed subrepresentation of $G$ consisting of bounded continuous functions $f \colon G \times M_0 \to k$ satisfying the following:
\begin{itemize}
\item[(i)] $f(g,am) = af(g,m)$ for any $(g,a,m) \in G \times k^{\circ} \times M_0$.
\item[(ii)] $f(g,m_1 + m_2) = f(g,m_1) + f(g,m_2)$ for any $(g,m_1,m_2) \in G \times M_0 \times M_0$.
\item[(iii)] $f(hg,m) = f(g,h^{-1}m)$ for any $(h,g,m) \in H \times G \times M_0$.
\end{itemize}
The inclusion $\m{Ind}_{H}^{G}(M_0^D) \hookrightarrow \m{C}_{\m{bd}}(G \times M_0,k)$ induces a continuous surjective $k^{\circ}[G]$-linear homomorphism $\Pi \colon \m{C}_{\m{bd}}(G \times M_0,k)^D \twoheadrightarrow \m{Ind}_{H}^{G}(M_0^D)^D$ by Hahn--Banach theorem (\cite{ST} Proposition 9.2). Beware that $\m{C}_{\m{bd}}(G,M_0^D)^D$ can not be expressed as a set of $M_0^D$-valued measures when $G$ is not compact and $M_0^D$ is not finite dimensional. On the other hand, since $k^{\circ}$ is compact, $\m{C}_{\m{bd}}(G \times M_0,k)^D$ coincides with the $k^{\circ}$-module $\m{M}(G \times M_0,k^{\circ})$ of integral $k$-valued measures on $G \times M_0$ endowed with the topology of pointwise convergence. Since the domain and the codomain of $\Pi$ are compact and Hausdorff, the codomain is homeomorphic to the coimage. Thus it suffices to determine $\ker \Pi$. For each $(g,m) \in G \times M$, we denote by $\delta_{g,m} \in \m{M}(G \times M_0,k^{\circ})$ the Dirac operator centred at $(g,m)$. We put $\mu^{M_0,\m{m}}_{g,a,m} \coloneqq a \delta_{g,m} - \delta_{g,am}$ for each $(g,a,m) \in G \times k^{\circ} \times M$, $\mu^{M_0,\m{a}}_{g,m,m'} \coloneqq \delta_{g,m+m'} - \delta_{g,m} - \delta_{g,m'}$ for each $(g,m,m') \in G \times M \times M$, and $\mu^{H,\rho}_{g,h,m} \coloneqq \delta_{hg,m} - \delta_{g,h^{-1}m}$ for each $(g,h,m) \in G \times H \times M$. We denote by $\mu^{M_0,\m{m}} + \mu^{M_0,\m{s}} + \mu^{H,\rho} \subset \m{M}(G \times M_0,k^{\circ})$ the closed $k^{\circ}$-submodule generated by the images of $\mu^{M_0,\m{m}}$, $\mu^{M_0,\m{s}}$, and $\mu^{H,\rho}$. We set
\begin{eqnarray*}
  \m{Ind}_{k^{\circ}[[H]]}^{k^{\circ}[[G]]}(M_0) \coloneqq \m{M}(G \times M_0,k^{\circ})/\mu^{M_0,\m{m}} + \mu^{M_0,\m{s}} + \mu^{H,\pi}
\end{eqnarray*}
and endow it with the quotient topology. We interpret the induction of a unitarisable Banach representation as the induction of a compact Hausdorff $k^{\circ}$-flat $k^{\circ}$-linear topological Iwasawa module.

\begin{thm}
\label{dual of the induction}
Let $V_0$ be a unitary Banach $k$-linear representation of $H$, and $M_0$ a compact Hausdorff $k^{\circ}$-flat $k^{\circ}$-linear topological $k^{\circ}[[H]]$-module with $V_0^D = M_0$. Then there is a natural $G$-equivariant isomorphism
\begin{eqnarray*}
  \m{Ind}_{H}^{G}(V_0)^D \cong \m{Ind}_{k^{\circ}[[H]]}^{k^{\circ}[[G]]}(M_0)
\end{eqnarray*}
of topological $k^{\circ}$-modules.
\end{thm}

In particular, the $k^{\circ}[G]$-module structure of the right hand side gives a topological $k^{\circ}[[G]]$-module structure.

\begin{proof}
It is easy to see that $\mu^{M_0,\m{m}} + \mu^{M_0,\m{s}} + \mu^{H,\pi} \subset \m{M}(G \times M_0,k^{\circ})$ is contained in the kernel of the canonical projection $\Pi \colon \m{M}(G \times M_0,k^{\circ}) \twoheadrightarrow \m{Ind}_{H}^{G}(V_0)^D$. Let $\mu \in \ker \Pi$. For any $f \in \m{C}_{\m{bd}}(G \times M_0,k)$ and any $\epsilon > 0$, we show that there is a $\mu' \in \mu^{M_0,\m{m}} + \mu^{M_0,\m{s}} + \mu^{H,\pi}$ such that $\v{\mu(f) - \mu'(f)} < \epsilon$. If $f \in \m{Ind}_{H}^{G}(V_0)$, then $\mu(f) = 0$, and hence $\mu' \coloneqq 0$ satisfies the desired inequality. If $f \notin \m{Ind}_{H}^{G}(V_0)$, then $f$ does not satisfies at least one of the conditions (i)-(iii) above. First, suppose that $f$ does not satisfy (i), and take a $(g,a,m) \in G \times k^{\circ} \times M_0$ such that $\alpha \coloneqq f(g,am) - af(g,m) \neq 0$. Set $\mu' \coloneqq \alpha^{-1} \mu(f) \mu^{M_0,\m{m}}_{g,a,m}$. Then $\mu'(f) = \mu(f)$, and hence $\v{\mu(f) - \mu'(f)} = 0 < \epsilon$. Next, suppose that $f$ does not satisfy (ii), and take a $(g,m,m') \in G \times M_0 \times M_0$ such that $\alpha \coloneqq f(g,m+m') - f(g,m) - f(g,m') \neq 0$. Set $\mu' \coloneqq \alpha^{-1} \mu(f) \mu^{M_0,\m{a}}_{g,m,m'}$. Then $\mu'(f) = \mu(f)$, and hence $\v{\mu(f) - \mu'(f)} = 0 < \epsilon$. Finally, suppose that $f$ does not satisfy (iii), and take a $(g,h,m) \in G \times H \times M_0$ such that $\alpha \coloneqq f(hg,m) - f(g,h^{-1}m) \neq 0$. Set $\mu' \coloneqq \alpha^{-1} \mu(f) \mu^{H,\rho}_{g,h,m}$. Then $\mu'(f) = \mu(f)$, and hence $\v{\mu(f) - \mu'(f)} = 0 < \epsilon$. Thus $\mu \in \mu^{M_0,\m{m}} + \mu^{M_0,\m{s}} + \mu^{H,\pi}$. We conclude that $\ker \Pi = \mu^{M_0,\m{m}} + \mu^{M_0,\m{s}} + \mu^{H,\pi}$, and hence $\Pi$ induces an isomorphism $\m{Ind}_{k^{\circ}[[H]]}^{k^{\circ}[[G]]}(M_0) \cong \m{Ind}_{H}^{G}(V_0)^D$.
\end{proof}

The functors $\m{Ind}_{H}^{G}$ and $\m{Ind}_{k^{\circ}[[H]]}^{k^{\circ}[[G]]}$ are $k^{\circ}$-linear, and so is the natural equivalence $D \circ \m{Ind}_{H}^{G} \circ D \cong \m{Ind}_{k^{\circ}[[H]]}^{k^{\circ}[[G]]}$ in Theorem \ref{dual of the induction}. Thus localising it, we obtain a $k$-linear natural equivalence $\mathbb{D} \circ \m{Ind}_{H}^{G} \circ D \cong (\m{Ind}_{k^{\circ}[[H]]}^{k^{\circ}[[G]]}) \otimes_{k^{\circ}} k$.

\begin{crl}
Let $V_0$ be a unitarisable Banach $k$-linear representation of $H$, and $M_0$ a compact Hausdorff $k^{\circ}$-flat $k^{\circ}$-linear topological $k^{\circ}[[H]]$-module with $V_0^{\mathbb{D}} \cong (M_0)_k$. Then there is a natural isomorphism
\begin{eqnarray*}
  \m{Ind}_{H}^{G}(V_0)^{\mathbb{D}} \cong \left( \m{Ind}_{k^{\circ}[[H]]}^{k^{\circ}[[G]]}(M_0) \right)_k
\end{eqnarray*}
of integral compactly generated locally convex $k[[G]]$-modules.
\end{crl}

\subsection{Continuous Parabolic Induction}
\label{Continuous Parabolic Induction}

We have constructed the dual notion of the continuous induction of a unitarisable Banach representation as an explicit quotient of the space of integral measures. We apply this result to a continuous parabolic induction on $\m{GL}_n(k)$.

\vspace{0.1in}
Let $H \subset G$ be a closed subgroup, and $K \subset G$ a compact subspace. Suppose $HK = G$. Then $H \backslash G$ is compact with respect to the quotient topology. For an object $M_0$ of $\m{TMod}_{\m{flat,lin}}^{\m{cpt,sep}}(k^{\circ}[[H]])$, we give a more practical expression of $\m{Ind}_{k^{\circ}[[H]]}^{k^{\circ}[[G]]}(M_0)$. Since $K \times M_0$ is compact, $\m{C}_{\m{bd}}(K \times M_0,k)$ coincides with the space $\m{C}(K \times M_0,k)$ of $k$-valued continuous functions on $K \times M$, and its dual coincides with $k^{\circ}[[K \times M_0]]$ endowed with the inverse limit topology $\tau_{\m{w}} = \tau_{\m{s}}$. We consider the composite $r_K \colon \m{Ind}_{H}^{G}(M_0^D) \to \m{C}(K \times M_0,k)$ of the embedding $\m{Ind}_{H}^{G}(M_0^D) \hookrightarrow \m{C}_{\m{bd}}(G \times M_0,k)$ and the restriction map $\m{C}_{\m{bd}}(G \times M_0,k) \to \m{C}(K \times M_0,k) \colon f \mapsto f|_{K \times M_0}$. Then $r_K$ is injective by the conditions (iii) and $HK = G$. The quotient norm on the coimage coincides with the norm restricted to the image because $H$ acts isometrically on $M_0^D$. Therefore we identify $\m{Ind}_{H}^{G}(M_0^D)$ as a closed $k$-linear subspace of $\m{C}(K \times M_0,k)$. We remark that if we start from a unitarisable Banach $k$-linear representation $V_0$ of $H$ instead of $M_0$, then $r_K$ is admissible in the sense that the quotient norm on the coimage is equivalent to the norm restricted to the image. The dual $\Pi_K \colon k^{\circ}[[K \times M_0]] \twoheadrightarrow \m{Ind}_{k^{\circ}[[H]]}^{k^{\circ}[[G]]}(M_0)$ of $r_k$ is surjective by Hahn--Banach theorem (\cite{ST} Proposition 9.2), and since $k^{\circ}[[K \times M_0]]$ and $\m{Ind}_{k^{\circ}[[H]]}^{k^{\circ}[[G]]}(M_0)$ are compact and Hausdorff, $\Pi_K$ is a quotient map.

\vspace{0.1in}
In addition, suppose that the multiplication $H \times K \to G$ is bijective. Then the inverse $p_1 \times p_2 \colon G \to H \times K$ is also continuous. Indeed, let $(g_i)_{i \in I}$ be a convergent net in $G$ with a limit $g \in G$. Since $K$ is compact, there is a subnet $(p_2(g_{i_j}))_{j \in J}$ of $(p_2(g_i))_{i \in I}$ with a limit $k \in K$. By the continuity of the multiplication and the inverse of $G$, the net $(p_1(g_{i_j}))_{j \in J} = (g_{i_j} p_2(g_{i_j})^{-1})_{j \in J}$ in $H$ converges to $gk^{-1} \in G$. Since $H \subset G$ is closed, $gk^{-1} \in H$. It follows from the bijectivity of $p_1 \times p_2$ that $k = p_2(g)$ and $gk^{-1} = p_1(g)$. This implies the continuity of $p_1 \times p_2$. Through the identification $G \cong H \times K$ as topological spaces, every $k$-valued continuous function on $K \times M_0$ is uniquely extended to a $k$-valued bounded continuous function on $G \times M_0$ satisfying the condition (iii). Therefore the image of $\m{Ind}_{H}^{G}(M_0^D) \hookrightarrow \m{C}(K \times M_0,k)$ coincides with functions satisfying the restriction of the conditions (i) and (ii). Moreover, the embedding is $G$-equivariant with respect to the isometric action $G \times \m{C}(K \times M_0,k) \to \m{C}(K \times M_0,k)$ given by setting $(gf)(k,m) \coloneqq f(p_2(kg),p_1(kg)^{-1}m)$ for each $(g,f,k,m) \in G \times \m{C}(K \times M_0,k) \times K \times M_0$. Since a continuous function on a compact uniform space is uniformly continuous, the action is strongly continuous. We endow $\m{M}(K \times M_0,k^{\circ})$ with the induced structure of a $k^{\circ}[[G]]$-module. We denote by $k^{\circ}[[K]] \hat{\otimes}_{k^{\circ}} M_0$ the completion of $k^{\circ}[[K]] \otimes_{k^{\circ}} M_0$ with respect to the linear topology generated by $k^{\circ}$-submodules of the form
\begin{eqnarray*}
  (L \otimes_{k^{\circ}} M_0) + (k^{\circ}[[K]] \otimes_{k^{\circ}} L') \subset k^{\circ}[[K]] \otimes_{k^{\circ}} M_0
\end{eqnarray*}
for open $k^{\circ}$-submodules $L \subset k^{\circ}[[K]]$ and $L' \subset M_0$. For a set $I$, we put $k^{\circ}[I] \coloneqq (k^{\circ})^{\oplus I}$ following the convention of the underlying $k^{\circ}$-module of a group algebra. We consider the action
\begin{eqnarray*}
  G \times (k^{\circ}[K] \otimes_{k^{\circ}} M_0) & \to & k^{\circ}[K] \otimes_{k^{\circ}} M_0 \\
  (g,[k] \otimes m) & \mapsto & [p_2(kg^{-1})] \otimes p_1(kg^{-1})^{-1}m
\end{eqnarray*}
of the underlying set of $G$. It is an action of a group because the isomorphism
\begin{eqnarray*}
  & & k^{\circ}[K] \otimes_{k^{\circ}} M_0 \cong k^{\circ}[K] \otimes_{k^{\circ}} k^{\circ}[H] \otimes_{k^{\circ}[H]} M_0 \cong (k^{\circ}[H] \otimes_{k^{\circ}} k^{\circ}[K]) \otimes_{k^{\circ}[H]} M_0 \\
  & \cong & k^{\circ}[H \times K] \otimes_{k^{\circ}[H]} M_0 \cong k^{\circ}[G] \otimes_{k^{\circ}[H]} M_0 
\end{eqnarray*}
of $k^{\circ}$-modules is $G$-equivariant with respect to the action on the right hand side as the induced representation. By the continuity of the multiplication $G \times G \to G$, $p_1 \times p_2$, and the action $G \times M_0 \to M_0$, the action on $k^{\circ}[K] \otimes_{k^{\circ}} M_0$ is strongly continuous with respect to the relative topology of $k^{\circ}[[K]] \hat{\otimes}_{k^{\circ}} M_0$, and hence it extends to a strongly continuous action
\begin{eqnarray*}
  G \times (k^{\circ}[[K]] \hat{\otimes}_{k^{\circ}} M_0) \to k^{\circ}[[K]] \hat{\otimes}_{k^{\circ}} M_0.
\end{eqnarray*}
Through the action, we endow $k^{\circ}[[K]] \hat{\otimes}_{k^{\circ}} M_0$ with a structure of topological $k^{\circ}[[G]]$-module by Proposition \ref{int iwa - rep cpt 3}.

\begin{prp}
In the situation above, there is a natural isomorphism
\begin{eqnarray*}
  \m{Ind}_{k^{\circ}[[H]]}^{k^{\circ}[[G]]}(M_0) \cong k^{\circ}[[K]] \hat{\otimes}_{k^{\circ}} M_0
\end{eqnarray*}
of topological $k^{\circ}[[G]]$-modules.
\end{prp}

\begin{proof}
By a similar argument as in the proof of Theorem \ref{dual of the induction}, $\ker \Pi_K$ coincides with the closed $k^{\circ}$-submodule generated by the images of $\mu^{M_0,\m{m}}$ and $\mu^{M_0,\m{a}}$ restricted to $K$. Since $K$ and $M_0$ are compact, any clopen covering of $K \times M_0$ admits a refinement by the direct product of the clopen coverings of $K$ and $M_0$. It implies that there is a natural isomorphism
\begin{eqnarray*}
  & & k^{\circ}[[K \times M_0]] \cong \varprojlim_{P,L'} \m{C}_0(P \times (M_0/L'),k^{\circ}) = \varprojlim_{P,L'} (k^{\circ})^{P \times (M_0/L')} \\
  & \cong & \varprojlim_{P,L'} \left( (k^{\circ})^P \otimes_{k^{\circ}} k^{\circ}[M_0/L'] \right) \cong \varprojlim_{L,L'} \left( (k^{\circ}[[K]]/L) \otimes_{k^{\circ}} k^{\circ}[M_0/L'] \right)
\end{eqnarray*}
of topological $k^{\circ}$-modules, where $P$, $L$, and $L'$ in the limit run through all $P \in \P(K)$, all open $k^{\circ}$-submodules $L \subset k^{\circ}[[K]]$, and all open $k^{\circ}$-submodules $L' \subset M_0$. We follows the same convention of $P$, $L$, and $L'$ as above. The last isomorphism holds because $(k^{\circ})^P \otimes_{k^{\circ}} k^{\circ}[M_0/L']$ is a finitely generated $k^{\circ}$-module. We have
\begin{eqnarray*}
  (k^{\circ}[[K]] \otimes_{k^{\circ}} M_0)/((L \otimes_{k^{\circ}} M_0) + (k^{\circ}[[K]] \otimes_{k^{\circ}} L')) \cong (k^{\circ}[[K]]/L) \otimes_{k^{\circ}} (M_0/L')
\end{eqnarray*}
by the flatness of the $k^{\circ}$-module $L$ and the right exactness of $(k^{\circ}[[K]]/L) \otimes_{k^{\circ}} (\cdot)$, and hence $k^{\circ}[[K]] \hat{\otimes}_{k^{\circ}} M_0$ is profinite. The kernel of the $k^{\circ}$-linear homomorphism
\begin{eqnarray*}
  (k^{\circ}[[K]]/L) \otimes_{k^{\circ}} k^{\circ}[M_0/L'] \to (k^{\circ}[[K]]/L) \otimes_{k^{\circ}} (M_0/L')
\end{eqnarray*}
induced by the scalar multiplication $k^{\circ}[M_0/L'] \to M_0/L'$ coincides with the image $(\ker \Pi_K)_{L,L'}$ of $\ker \Pi_K$ by the right exactness of $(k^{\circ}[[K]]/L) \otimes_{k^{\circ}} (\cdot)$. Therefore taking the inverse limit of the exact sequence
\begin{eqnarray*}
  0 \to (\ker \Pi_K)_{L,L'} \to (k^{\circ}[[K]]/L) \otimes_{k^{\circ}} k^{\circ}[M_0/L'] \to (k^{\circ}[[K]]/L) \otimes_{k^{\circ}} (M_0/L') \to 0
\end{eqnarray*}
of finitely generated $k^{\circ}$-modules, we obtain an exact sequence
\begin{eqnarray*}
  0 \to \ker \Pi_K \to k^{\circ}[[K \times M_0]] \to k^{\circ}[[K]] \hat{\otimes}_{k^{\circ}} M_0
\end{eqnarray*}
of profinite $k^{\circ}$-modules by the closedness of $\ker \Pi_K \subset k^{\circ}[[K \times M_0]]$. Since $\ker \Pi_K$ is profinite, the system $((\ker \Pi_K)_{L,L'})$ satisfies the Mittag-Leffler condition, and hence the sequence
\begin{eqnarray*}
  0 \to \ker \Pi_K \to k^{\circ}[[K \times M_0]] \to k^{\circ}[[K]] \hat{\otimes}_{k^{\circ}} M_0 \to 0
\end{eqnarray*}
of $k^{\circ}$-modules is exact. Each homomorphism is continuous, and since $k^{\circ}[[K \times M_0]]$ and $k^{\circ}[[K]] \hat{\otimes}_{k^{\circ}} M_0$ are compact and Hausdorff, it gives an isomorphism
\begin{eqnarray*}
  \m{Ind}_{k^{\circ}[[H]]}^{k^{\circ}[[G]]}(M_0) \cong k^{\circ}[[K \times M_0]]/\ker \Pi_K \cong k^{\circ}[[K]] \hat{\otimes}_{k^{\circ}} M_0
\end{eqnarray*}
of topological $k^{\circ}$-modules. It is $G$-equivariant by the definition of the action on the right hand side, and is an isomorphism of topological $k^{\circ}[[G]]$-modules.
\end{proof}

\begin{exm}
For an $n \in \N$, denote by $B_n^{+}(k) \subset \m{GL}_n(k)$ the Borel subgroup of upper triangular matrices, by $K_n^{-} \subset \m{GL}_n(k)$ the compact subset of the normalisations of lower triangular matrices with diagonal entries $1$, and by $\mathscr{S}_n \subset \m{GL}_n(k)$ is the finite subgroup of permutations of the canonical basis. Here the normalisation of a lower triangular matrix means the unique lower triangular matrix in $\m{M}_n(k^{\circ}) \cap \m{GL}_n(k)$ obtained by multiplying an invertible diagonal matrix from left such that for each row, an entry with norm $1$ exists and the entry in the greatest column among entries with norm $1$ is $1$. By a method of LUP-decomposition, $\m{GL}_n(k)$ is expressed as the product $B_n^{+}(k) K_n^{-} \mathscr{S}_n$, and the multiplication $B_n^{+}(k) \times K_n^{-} \mathscr{S}_n \to \m{GL}_n(k)$ is bijective. Thus for a compact Hausdorff $k^{\circ}$-flat $k^{\circ}$-linear topological $k^{\circ}[[B_n^{+}(k)]]$-module $M_0$, we have a natural isomorphism
\begin{eqnarray*}
  \m{Ind}_{k^{\circ}[[B_n^{+}(k)]]}^{k^{\circ}[[\m{GL}_n(k)]]}(M_0) \cong k^{\circ}[[K_n^{-} \mathscr{S}_n]] \hat{\otimes}_{k^{\circ}} M_0 \cong (k^{\circ}[[K_n^{-}]] \hat{\otimes}_{k^{\circ}} M_0)^{\mathscr{S}_n}
\end{eqnarray*}
of topological $k^{\circ}[[\m{GL}_n(k)]]$-modules. Beware that the $k^{\circ}[[\m{GL}_n(k)]]$-module structures of the right hand sides depend on the decomposition $\m{GL}_n(k) = B_n^{+}(k) K_n^{-} \mathscr{S}_n$. Applying this result to a continuous parabolic induction, for a unitary Banach $k$-linear representation $V_0$ of $B_n^{+}(k)$, we have a natural isometric isomorphism
\begin{eqnarray*}
  \m{Ind}_{B_n^{+}(k)}^{\m{GL}_n(k)}(V_0) \cong \left( k^{\circ}[[K_n^{-} \mathscr{S}_n]] \hat{\otimes}_{k^{\circ}} V_0^D \right)^D \cong \left( \left( k^{\circ}[[K_n^{-}]] \hat{\otimes}_{k^{\circ}} V_0^D \right)^D \right)^{\mathscr{S}_n}
\end{eqnarray*}
of unitary Banach $k$-linear representations of $\m{GL}_n(k)$.
\end{exm}

\vspace{0.4in}
\newpage
\addcontentsline{toc}{section}{Acknowledgements}
\noindent {\Large \bf Acknowledgements}
\vspace{0.1in}

I am extremely grateful to Takeshi Tsuji for constructive advices in seminars. I thank my friends for daily discussions. I greatly appreciate my family's deep affection.

\end{document}